\documentclass[12pt]{article}
\usepackage{amssymb}
\usepackage{amsmath}
\usepackage{amsthm}
\usepackage{amsfonts}
\usepackage{graphicx}
\usepackage{subcaption}
\usepackage{enumerate}
\usepackage{natbib}
\usepackage{color}
\usepackage{algorithm}
\usepackage{algpseudocode}
\usepackage{url} 
\usepackage[font=footnotesize]{caption} 
\usepackage{float}
\newcommand{\blind}{1}

\newcommand{\RNum}[1]{\uppercase\expandafter{\romannumeral #1\relax}}
\newcommand{\bx}{x}
\newcommand{\bz}{z}
\newcommand{\by}{{\bf y}}

\DeclareMathOperator{\Tr}{Tr}
\newcommand{\ve}{\varepsilon}

\newtheorem{Prop}{Proposition}
\newtheorem{theorem}{Theorem}

\begin{document}

\def\spacingset#1{\renewcommand{\baselinestretch}%
{#1}\small\normalsize} \spacingset{1}


\if1\blind
{
  \title{\bf Optimal Design of Experiments on Riemannian Manifolds}
  \author{\small Hang Li\thanks{
    The authors gratefully acknowledge \textit{NSF grant CMII 1537987}}\hspace{.2cm}\\
    \small Department of Industrial and Manufacturing Engineering, \\ 
    \small Pennsylvania State University\\
    \small and \\
    \small Enrique Del Castillo \\
    \small Department of Industrial and Manufacturing Engineering and Dept. of Statistics, \\ 
    \small Pennsylvania State University}
    \date{\small November 25, 2019}
  \maketitle
} \fi

\if0\blind
{
  \bigskip
  \bigskip
  \bigskip
  \begin{center}
    {\LARGE\bf Optimal Design of Experiments on Riemannian Manifolds}
\end{center}
  \medskip
} \fi

\begin{abstract}
The theory of optimal design of experiments has been traditionally developed on an Euclidean space. In this paper, new theoretical results and an algorithm for finding the optimal design of an experiment located on a Riemannian manifold are provided. It is shown that analogously to the results in Euclidean spaces, D-optimal  and G-optimal designs are equivalent on manifolds, and we provide a lower bound for the maximum prediction variance of the response evaluated over the manifold. In addition, a converging algorithm that finds the optimal experimental design on manifold data is proposed. Numerical experiments demonstrate the importance of considering the manifold structure in a designed experiment when present, and the superiority of the proposed algorithm.   
\end{abstract}

\noindent%

\spacingset{1.5} 
\section{Introduction}\label{sec:intro}

Supervised learning models typically need to be trained on large amounts of labeled instances to perform well. While many modern systems can easily produce a large number of unlabeled instances at low cost, the labeling process can be very difficult, expensive or time-consuming. For example, audio data require experienced linguists to spend much longer time than the audio itself to precisely annotate the speech utterances. Given a learning model, nonidentical labeled instances contain different amounts of information and contribute to the learning process in different ways. Therefore, an interesting and practical question arises: how to choose the most informative instances to label so that one can improve the learning rate of the model and reduce the labeling cost at the same time?

In statistics, the learning problem of selecting which instances to label is closely related to classical Design of Experiments   \citep[DOE,][]{Wu2009}. Traditional DOE was developed for physical experiments in agricultural applications, where the goal is to explore the relationship between several input covariates and one output response under limited experimental resources. An important question widely studied within DOE is to define optimality criteria for experimental design and algorithms to obtain such designs, that is, Optimal Design of Experiments ({ODOE}). The goal of ODOE is to develop experimental designs that are optimal with respect to some statistical criterion. In the classical theory of ODOE, a linear model is usually assumed and the statistical criteria are typically related to the model parameter estimates or to the model predictions  \citep[see][]{kiefer_wolfowitz_1960, Fedorov1972,Pukelsheim2006,FedorovLeonov2013}. 
In traditional DOE problems, the number of covariates or  ``factors" of interest in an experiment is relatively small and the experimental region is usually assumed to be Euclidean. However, in some modern  learning tasks, such as image recognition and text categorization, the dimension of text or image data is often much higher than the dimension of covariates in a traditional agricultural or industrial experiment. 

In order to perform a statistical learning task under these conditions, a manifold hypothesis is made, which assumes that although  
the training data of interest are available in a high-dimensional ambient space, there exists a lower-dimensional manifold where the data are located. This is in contrast to traditional linear dimensionality reduction based on principal components, 
where a linear subspace is assumed with the hope the data is concentrated in it.
 The manifold hypothesis is often observable in  high-dimensional data. Starting with the the work by  \cite{RS00S} and  \cite{TSL00S}, a wide body of literature has shown how high-dimensional data, such as text or image data, frequently lie on a lower dimensional manifold. 
 and is usually sparse in its high-dimensional ambient space  \citep[see, e.g.,][]{CWJASA2013,LTZD2017JASA,nature2018,doi:10.1080/01621459.2019.1610660}. 
From a DOE point of view, it is infeasible to obtain and label enough training instances to fill up the high-dimensional ambient space. However, one could select points from a lower dimensional manifold space if the data points were much more dense on this low-dimensional space. Unfortunately,  traditional DOE methods fail to take into account these complex characteristics of modern high-dimensional data \citep{LCR2019TEST}. 

\begin{figure}[H]
        \centering
        \begin{subfigure}[b]{0.31\textwidth}
            \centering
            \includegraphics[width=\textwidth]{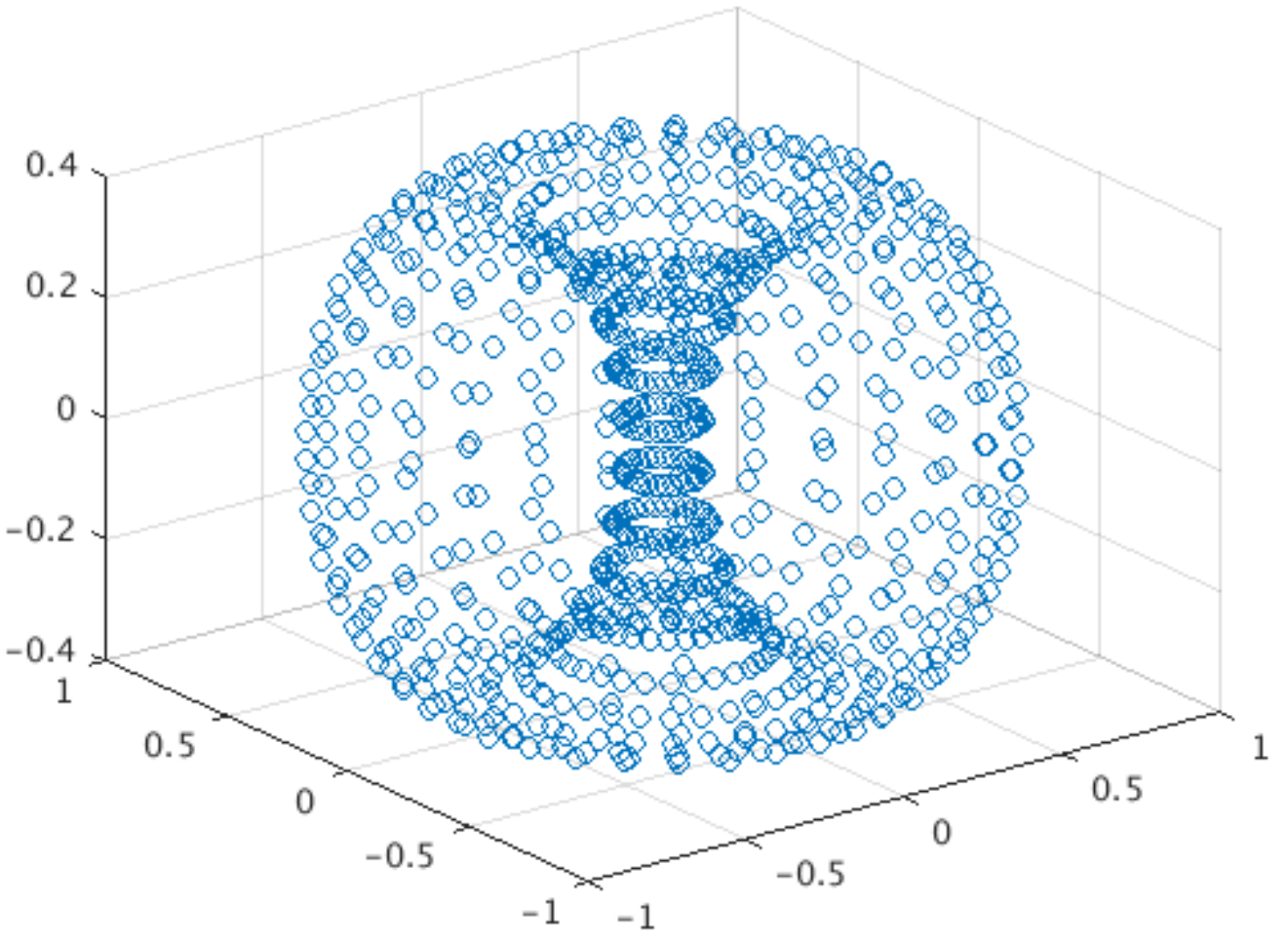}
            \caption{Torus Data}
        \end{subfigure}
        \quad
        \begin{subfigure}[b]{0.31\textwidth}
            \centering
            \includegraphics[width=\textwidth]{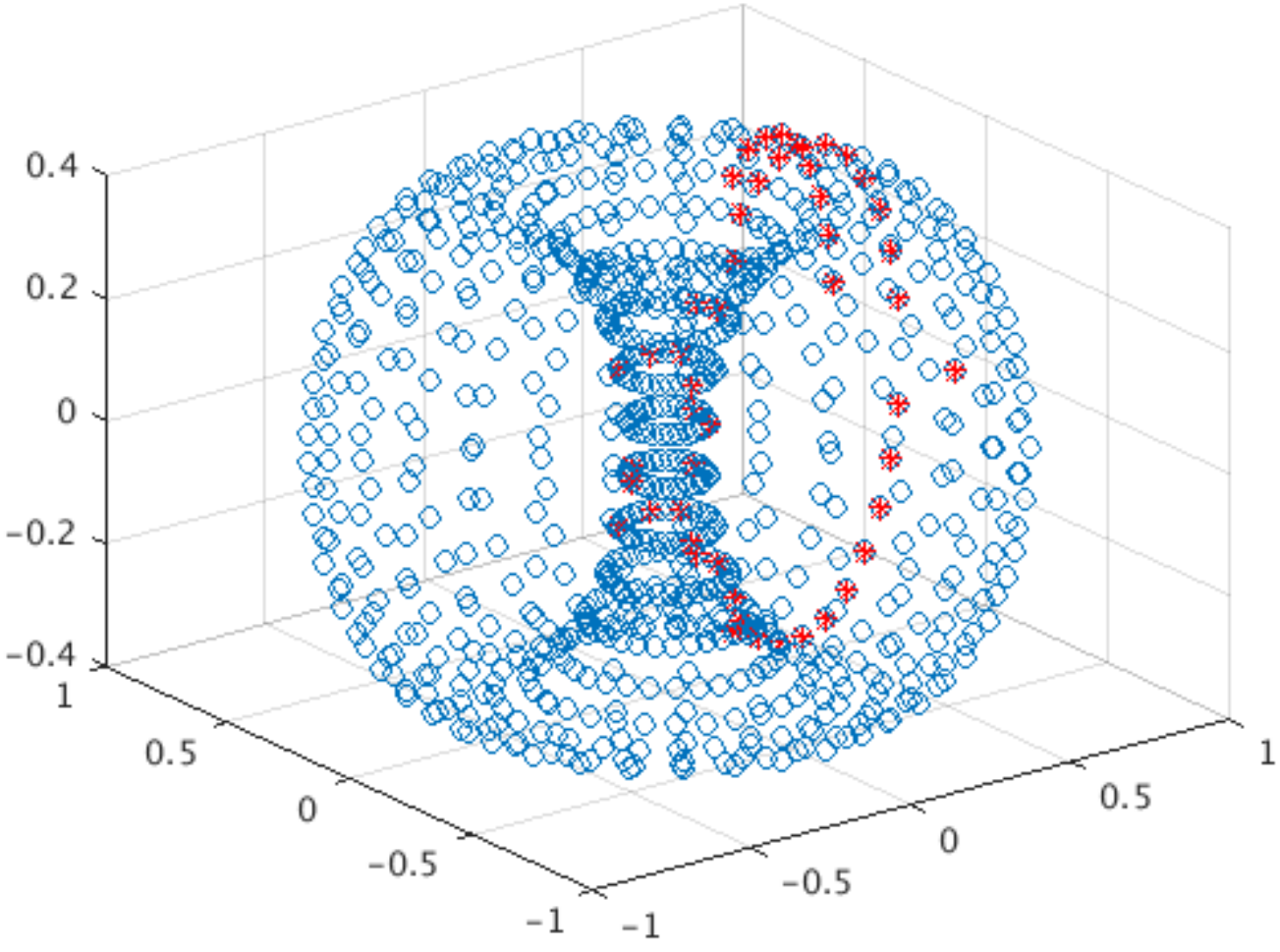}
            \caption{Experimental design \RNum{1}}
        \end{subfigure}
        \quad
        \begin{subfigure}[b]{0.31\textwidth}
            \centering
            \includegraphics[width=\textwidth]{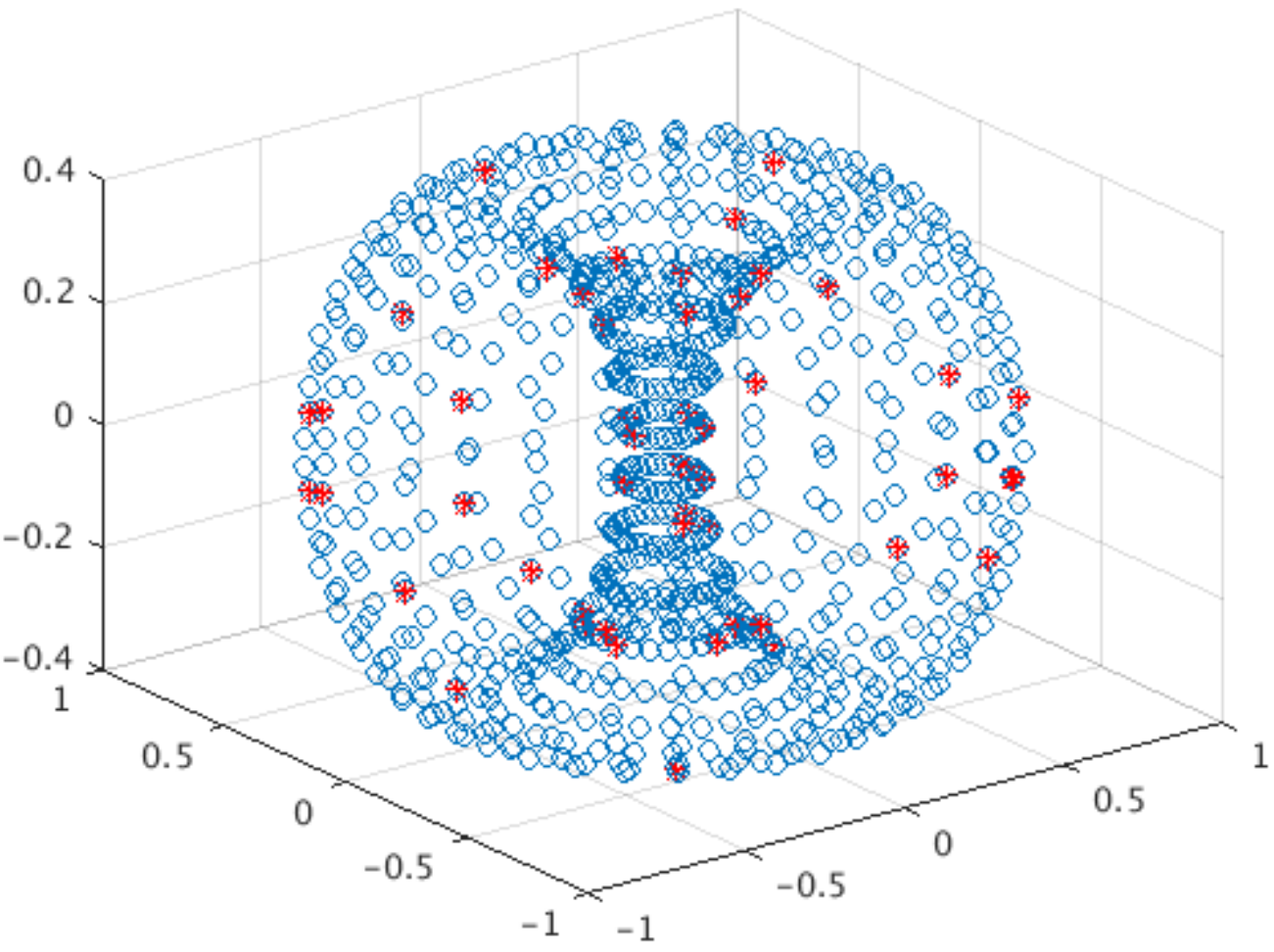}
            \caption{Experimental design \RNum{2}}
        \end{subfigure}
        \caption{A toy example that data lie on a 2-dimensional Torus embedded in a 3-dimensional ambient space. The red points are the instances selected to be labeled in the corresponding experiments.}
        \label{fig:intro}
\end{figure}

Figure \ref{fig:intro}(a) demonstrates an easy to visualize example, where the data points are available in a 3-dimensional Euclidean space but truly lie on the 2-dimensional surface of a Torus. Two different experimental designs on this dataset are provided in Figure \ref{fig:intro}(b) and \ref{fig:intro}(c). Different designs will lead to different learning performance. The motivating question is therefore: how to find the optimal design that improves the learning performance the most, while incorporating the  manifold structure where data lie into account? 

The goal of this paper is to develop theory and an algorithm for constructing optimal experimental designs on high-dimensional manifold data, which minimize the number of experimental runs  and at the same time acquire as much useful information about the response as possible. We assume training data are located on a lower dimensional Riemannian manifold, loosely defined as a curved space which when seen over a sufficiently small neighborhood resembles Euclidean (flat) space. 

Although some previous authors have implemented ODOE criteria as Active Learning strategies for high-dimensional data \citep{H10IEEETIP,CCBWZZ2010,Alaeddini2019IISE}, as far as we know, no existing work has provided theoretical guarantees of such  experimental designs on Riemannian manifolds. 

Our contributions are summarized as follows: a) we prove a new Equivalence Theorem for continuous optimal designs on manifold data, which shows how a D-optimal and a G-optimal designs are equivalent on Riemannian manifolds; b) we provide a new lower bound for the maximum prediction variance over the manifold and show how the lower bound can be achieved by a D/G optimal design; c) we propose a new algorithm, ODOEM (Optimal Design of Experiments on Manifolds), for finding a continuous D/G optimal design on a Riemannian manifold, and prove that it is guaranteed to converge to the global D/G optimal design, and finally, d) we illustrate the superior performance of our ODOEM algorithm on both of synthetic manifold datasets and a real-world image dataset.

The rest of this paper is organized as follows. In Section \ref{sec:ODOE}, we briefly review the traditional ODOE problem on Euclidean space, and then introduce the \emph{manifold regularization} model of \cite{BNS06JMLR} on which our results are based, explaining the ODOE problem on manifolds. Section \ref{sec:theory} provides the theoretical justification behind our ODOEM algorithm, where a new equivalence theorem is given for Riemannian manifolds.  Section \ref{sec:alg} gives the proposed ODOEM algorithm and provides a convergence analysis. Finally, section \ref{sec:numerical} presents several numerical experiments conducted to demonstrate the effectiveness of the proposed algorithm for finding optimal designs on manifold data. We conclude the paper with a summary and some possible further research directions in section \ref{sec:conc}.


\section{Optimal Design of Experiments on Manifolds}\label{sec:ODOE}
\subsection{Traditional ODOE on Euclidean Space}\label{}
Consider initially a linear regression model
\begin{equation}\label{eqn:lr}
y = f(\bx,\beta) + \ve\;  = \; \beta^\top g(\bx) + \ve, 
\end{equation}
where $g: \mathbb{R}^d \rightarrow \mathbb{R}^p$ is some nonlinear function that maps from  the input space $\bx \in \mathbb{R}^d$ to the feature space $\mathbb{R}^p$, $\beta \in \mathbb{R}^p$ is a column vector of unknown parameters, and $\ve$ is assumed to have a $N(0,\sigma^2)$ distribution . Given a sample of $n$ design points $\{\bx_i\}_{i=1}^n$, if the corresponding response values $\{y_i\}_{i=1}^n$ are available, the well-known ordinary least squares estimates of the $\beta$ parameters are given by:
\begin{equation}\label{OLS_est}
    \hat \beta= \operatornamewithlimits{argmin}_{\beta\in \mathbb{R}^p} \left\{ \sum_{i=1}^n ( y_i - \beta^\top g(\bx_i))^2  \right\} \; = \; (X^\top X)^{-1} X^\top Y
\end{equation}
where $X$ is a $n \times p$ design matrix with $i$-th row defined as $g(\bx_i)^\top$, and $Y$ is a $n \times 1$ response vector. As a result, the corresponding fitted function is $\hat f(\bx) = \hat \beta^\top g(\bx)$.

Classical work on ODOE was developed by \cite{kiefer_wolfowitz_1960} and summarized by \cite{Fedorov1972} \citep[see also][]{Pukelsheim2006,FedorovLeonov2013}. Examples of optimality criteria for the linear regression model (\ref{eqn:lr}) are the \textit{D-optimality} criterion which minimizes the determinant of the covariance matrix of the parameter estimates $\operatorname{Var}(\hat \beta)= \sigma^2 (X^\top X)^{-1}$, and the \textit{G-optimality} criterion which minimizes the maximum prediction variance $\operatornamewithlimits{max}_{i=1,...,n}\big\{\operatorname{Var}(\hat{y}_i) \big\}$. These and similar criteria are called  ``\textit{alphabetic optimality}" design criteria by \cite{box2007response}.

While there have been recent attempts at applying alphabetic optimality criteria to manifold learning models \citep{H10IEEETIP,CCBWZZ2010,Alaeddini2019IISE}, no theoretical justification exists, as far as we know, to these methods, and no guarantees can be given for their success other than empirical experimentation. A new  theory for optimal experimental design is therefore needed that explicitly considers  high-dimensional manifold data, justify existing methods if possible, and that provides a principled way to develop new algorithms. Before we discuss the design of experiments on manifolds, first we need to introduce a manifold learning model by \cite{BNS06JMLR} that will be used in the sequel.

\subsection{Manifold Regularization Model}

In the standard paradigm of machine learning, learning is understood as a process that uses the training data $\{\bx_i\}_{i=1}^{n}$ to construct a function $f: \mathcal{X} \rightarrow \mathbb{R}$ that maps a data instance $\bx$ to a label variable $y$. Let $P$ be the joint distribution that generates labeled data $\{(\bx_i,y_i)\}_{i=1}^{l} \subset \mathcal{X}\times \mathbb{R}$ and $P_\mathcal{X}$ be the marginal distribution that generates unlabeled data {$\{\bx_i\}_{i=l+1}^{n} \subset \mathcal{X} \subset \mathbb{R}^d$}. In order to extend the learning of functions to general Riemannian manifolds, \cite{BNS06JMLR} assume that the conditional distribution $P(y|\bx)$ varies smoothly as $\bx$ moves along a manifold that supports $P_\mathcal{X}$. In other words, if two data points $\bx_1, \bx_2 \in \mathcal{X}$ are close as measured by an intrinsic (or geodesic) distance on this manifold, then the two probabilities of the labels, $P(y|\bx_1)$ and $P(y|\bx_2)$, will be similar. These authors developed a semi-supervised learning framework that involves solving the following double regularized objective function:
\begin{equation}\label{eqn:obj1}
  \hat{f}=\operatornamewithlimits{argmin}\limits_{f\in \mathcal{H}_\mathcal{K}} \left\{ \sum_{i=1}^{l} V(\bx_i,y_i,f)+\lambda_A \|f\|_{\mathcal{H}_\mathcal{K}}^2+\lambda_I \|f\|_I^2  \right\}
\end{equation}
where $V$ is a given loss function (such as squared loss $(y_i-f(\bx_i))^2$), $\mathcal{H}_\mathcal{K}$ is a \emph{Reproducing Kernel Hilbert Space} \citep[RKHS,][]{A50TAMS} with associated Mercer kernel $\mathcal{K}$, $ \| f \|_{\mathcal{H}_\mathcal{K}}^2$ is a penalty term with the norm of $\mathcal{H}_\mathcal{K}$ that imposes smoothness conditions in the ambient space \citep{W90SM}, and $\|f\|_I^2$ is a penalty term for non-smoothness along geodesics on the intrinsic manifold structure of $P_\mathcal{X}$. Moreover, $\lambda_A$ and $\lambda_I$ are two regularization parameters that control the amount of penalization in the ambient space and in the intrinsic manifold that supports $P_\mathcal{X}$, respectively. Recent work on non-Euclidean data methods are related to (\ref{eqn:obj1}), for instance, the spatial regression model proposed by \cite{ESP2016Bio} can be seen as the manifold regularization model (\ref{eqn:obj1}) without the ambient space regularization. While there are also numerous nonparametric regression models on manifolds \citep[see, e.g.,][]{CWJASA2013,MPT2014JASA,LTZD2017JASA}, we focus on this paper on the manifold regularization model from \cite{BNS06JMLR} since it provides a nice representer theorem, an advantage that will be clear in what follows.

Intuitively, the choice of $\|f\|^2_I$ should be a smoothness penalty corresponding to the probability distribution $P_\mathcal{X}$. However, in most real-world applications $P_\mathcal{X}$ is not known, and therefore empirical estimates of the marginal distribution must be used. Considerable research has been devoted to the case when $P_\mathcal{X}$ is supported on a compact manifold $\mathcal{M} \subset \mathbb{R}^d$ \citep{RS00S,TSL00S,BN03NC,DG03PNAS,CLLMNWZ05PNAS}. Under this assumption, it can be shown \citep[see][]{B03UCTHESIS,L04YALETHESIS} that problem (\ref{eqn:obj1}) can be reduced to
\begin{eqnarray}\label{eqn:mrobj2}
\hat f = \operatornamewithlimits{argmin}\limits_{f\in \mathcal{H}_\mathcal{K}} \left\{ \sum_{i=1}^{l} V(\bx_i,y_i,f)+\lambda_A \|f\|_{\mathcal{H}_\mathcal{K}}^2+\lambda_I \mathbf{f}^\top L \mathbf{f} \right\}
\end{eqnarray}
where $\mathbf{f}=[f(\bx_1),...,f(\bx_{n})]^\top$ and $L$ is the Laplacian matrix associated with the data adjacency graph $\mathcal{G}$ that is constructed on all the labeled and the unlabeled data points $\{\bx_i\}_{i=1}^{n}$. In particular, the graph Laplacian $L$ approximates the \textit{Laplace-Beltrami} operator acting on the continuous Riemannian manifold $\mathcal{M}$ \citep[see][]{BN05COLT,CLLMNWZ05PNAS,HAL05COLT}. The convergence of the graph Laplacian provides a theoretical justification to the common practice in manifold learning of using a graph and the corresponding geodesic distances as an approximate representation of the manifold $\mathcal{M}$, providing a precise sense in which the graph approaches $\mathcal{M}$ as the number of data points gets denser. This way, the term $\mathbf{f}^\top L \mathbf{f}$ serves as an approximation for $\|f\|_I^2 $, and enforces the penalization on the lack of smoothness of $f$ as it varies between adjacent points in the graph $\mathcal{G}$. 

The solution of the infinite dimensional problem (\ref{eqn:mrobj2}) can be represented in terms of a finite sum over the labeled and unlabeled points:
\begin{eqnarray} \label{representer}
 f(\bx) = \sum_{i=1}^{n}  \alpha_i \mathcal{K}(\bx_i,\bx)
\end{eqnarray}
where $\mathcal{K}(\cdot,\cdot)$ is the Mercer kernel associated with the ambient space $\mathcal{H}_\mathcal{K}$. This constitutes a representer theorem for problem (\ref{eqn:mrobj2}), similar to that in the theory of splines \citep{kimeldorf1970,W90SM}. 

\subsection{Regularized ODOE on Manifolds}

\cite{Vuchkov1977} provided the first discussion of a regularized method in the literature on ODOE, based on the ridge regression estimator:
\begin{equation}\label{ridge}
\hat \beta_{\text{ridge}} = \operatornamewithlimits{argmin}_{\beta\in \mathbb{R}^p} \left\{ \sum_{i=1}^l ( y_i - \beta^\top g(\bx_i))^2 + \lambda_{\text{ridge}} \| \beta\|^2  \right\} 
\end{equation}
Vuchkov's motivation was to use the ridge estimator to solve the singular or ill-conditioned problems that exist in the sequential application of a D-optimal design algorithm when the number of design points is smaller than the number of parameters to estimate. The ridge solution (\ref{ridge}) can be seen as a particular case of the more general learning problem (\ref{eqn:mrobj2}) where $V$ is a squared-loss function, the RKHS $\mathcal{H}_K$ is equipped with a $L^2$-norm  and the manifold regularization parameter $\lambda_I$ is zero.

To discuss the optimal experimental design for the general manifold regularization model (\ref{eqn:mrobj2}), we first introduce some notation. Without loss of generality, assume a sequential experimental design problem, starting with no labeled data at the beginning of the sequence. Let $\{\bz_i\}_{i=1}^k \subset \{\bx_i\}_{i=1}^n$ be the set of points that has been labeled at the $k$-th iteration, and $\textbf{y}=(y_1,...,y_k)^\top$ be the corresponding vector of responses or labels. Given a square loss function, the manifold regularization model (\ref{eqn:mrobj2}) becomes the Laplacian Regularized Least Squares ({LapRLS}) problem \citep{BNS06JMLR}:
\begin{equation}\label{eqn:LapRLS}
\hat f = \operatornamewithlimits{argmin}\limits_{f\in \mathcal{H}_K} \left\{ \sum_{i=1}^{k}(y_i-f(\bz_i))^2+\lambda_A \|f\|_{\mathcal{H}_\mathcal{K}}^2+\lambda_I \textbf{f}^\top  L \textbf{f} \right\}.
\end{equation}
Substituting the representer theorem solution (\ref{representer}) into (\ref{eqn:LapRLS}), we get a convex differentiable objective function with respect to $\alpha$:
\begin{equation}\label{eqn:obj_alpha}
\hat{\alpha}  = \operatornamewithlimits{argmin}\limits_{\alpha\in \mathbb{R}^n} \left\{ (\by-K_{XZ}^\top \alpha)^\top (\by-K_{XZ}^\top \alpha) +\lambda_A \alpha^\top K \alpha +\lambda_I \alpha^\top KLK \alpha \right\},
\end{equation}
where $K_{XZ}$ and $K$ are the Gram matrices defined by
\begin{eqnarray*}
K_{XZ}=
\begin{bmatrix}
\mathcal{K}(\bx_1,\bz_1) &  ... & \mathcal{K}(\bx_1,\bz_k)\\
\vdots &  \ddots & \vdots \\
\mathcal{K}(\bx_n,\bz_1) &  ... & \mathcal{K}(\bx_n,\bz_k)
\end{bmatrix}_{n \times k}, \; 
K=
\begin{bmatrix}
\mathcal{K}(\bx_1,\bx_1) &  ... & \mathcal{K}(\bx_1,\bx_n)\\
\vdots &  \ddots & \vdots \\
\mathcal{K}(\bx_n,\bx_1) &  ... & \mathcal{K}(\bx_n,\bx_n)
\end{bmatrix}_{n \times n},
\end{eqnarray*}
and $\mathcal{K}$ is the kernel embedded in the RKHS $\mathcal{H}_\mathcal{K}$. Taking the derivative of (\ref{eqn:obj_alpha}) with respect to $\alpha$ and making it equal to 0, we arrive at the following expression:
\begin{eqnarray}
 \hat{\alpha} = (K_{XZ}K_{XZ}^\top +\lambda_A  K + \lambda_I KLK)^{-1} K_{XZ}\mathbf{y}
\end{eqnarray}

Consider a linear model of the form (\ref{eqn:lr}) and a linear kernel for $\mathcal{H}_K$, the regression parameters $\beta$ can be estimated by
\begin{eqnarray}\label{beta_est1}
 \hat{\beta}=X^\top \hat{\alpha} = X^\top (XZ_k^\top Z_k X^\top+\lambda_A XX^\top + \lambda_I XX^\top L XX^\top)^{-1}XZ^\top \mathbf{y}
\end{eqnarray}
where
\begin{eqnarray}
Z_k=\left [   \begin{array}{c}
                  g(\bz_1)^\top\\
                  \vdots \\
                  g(\bz_k)^\top
                \end{array}
              \right], \; X=\left [   \begin{array}{c}
                  g(\bx_1)^\top\\
                  \vdots \\
                  g(\bx_n)^\top
                \end{array}
              \right], \; \mathbf{y}=\left [   \begin{array}{c}
                  y_1\\
                  \vdots \\
                  y_k
                \end{array}
              \right].
\end{eqnarray}
By some simple linear algebra (a formal proof is provided in the Appendix), the estimated parameters $\hat{\beta}$ can be simplified to
\begin{eqnarray}\label{beta_est2}
 \hat{\beta}=(Z_k^\top Z_k + \lambda_A I_p + \lambda_I X^\top LX)^{-1} Z^\top \mathbf{y}
\end{eqnarray}
Similarly to the theory of ODOR on Euclidean space, the regularized estimator resembles a Bayesian linear estimator, with the difference being that the regularization comes from the manifold penalization of high-dimensional data instead of some a priori covariance estimate \citep{Pukelsheim2006}.

\cite{H10IEEETIP} demonstrated that the covariance matrix of (\ref{beta_est2}) can be approximated as:  
\begin{equation}\label{eqn:cov_appr} 
\mathrm{Cov} (\mathbf{\hat{\beta}}) \approx \sigma^2 (Z_k^\top Z_k +\lambda_A I_p + \lambda_I X^\top L X)^{-1}.
\end{equation}
The determinant of covariance matrix (\ref{eqn:cov_appr}) is the statistical criterion we will minimize to obtain a D-optimal design for manifold data. Before we discuss the optimal design algorithm,  first we will provide its main theoretical justification.
\section{Equivalence Theorem on Manifolds}\label{sec:theory}
When the determinant of $Z_k^\top Z_k + \lambda_A I_p + \lambda_I X^\top L X$ is maximized, one obtains a D-optimal experimental design. In Euclidean space, ODOE indicates an equivalence between the D-optimality criteria and the G-optimality criteria, which minimizes the maximum prediction variance, as stated by the celebrated Kiefer-Wolfowitz ( {KW}) theorem \citep{kiefer_wolfowitz_1960,Kiefer1974}. In analogy with the KW theorem, in this section we aim to develop a new equivalence result for optimal experimental design based on the manifold regularization model (\ref{eqn:mrobj2}), which can then be used to justify algorithms for designing an optimal experiment on a Riemannian manifold.

Assume there is an infinite number of points $x$ that are uniformly distributed on a Riemannian manifold $\mathcal{M}$. Let $\epsilon$ be a continuous design on $\mathcal{M}$. For any continuous design $\epsilon$, based on the Carath\'{e}odory Theorem, it is known \citep[see][]{Fedorov1972} that $\epsilon$ can be represented as
\begin{eqnarray}
\epsilon =\left \{   \begin{array}{cccc}
                  z_1, z_2, ..., z_{n_0}\\
                  q_1,q_2,...,q_{n_0}
                \end{array}
              \right\}, \; \mathrm{where} \; \sum_{i=1}^{n_0} q_i=1.
\end{eqnarray} For any $\epsilon$, the corresponding information matrix of LapRLS model is defined as
\begin{equation} \label{Mcont}
 M_{Lap}(\epsilon)=\int_{z \in \mathcal{X}} \xi(z)g(z)g(z)^\top dz + \lambda_A I_p + \lambda_I \int_{x \in \mathcal{M}} g(x) \Delta_{ \mathcal{M}} g(x)^\top d\mu,  
\end{equation}
where $\xi$ is a probability measure of design $\epsilon$ on the experimental region $\mathcal{X}\subseteq \mathcal{M} \subset \mathbb{R}^p$, $\Delta_{ \mathcal{M}}$ is the \textit{Laplace-Beltrami} operator on $\mathcal{M}$, and $\mu$ is the uniform measure on $\mathcal{M}$. Note that the last two terms in (\ref{Mcont}) are independent of the design $\epsilon$, thus for simplicity, define 
\begin{equation} 
C = \lambda_A I_p + \lambda_I \int_{x \in \mathcal{M}} g(x) \Delta_{ \mathcal{M}} g(x)^\top d\mu.
\end{equation}
Then (\ref{Mcont}) can be written as 
\begin{equation} 
  M_{Lap}(\epsilon)=\int_{z \in \mathcal{X}} \xi(z)g(z)g(z)^\top dz + C.  
\end{equation}
Based on the parameters estimates (\ref{beta_est2}), for a given continuous design $\epsilon$, the prediction variance at a test point $z$ is
\begin{eqnarray}
d(z,\epsilon)=\mathrm{Var}\Big[\mathbf{\hat{\beta}}^\top g(z) \Big]=g(z)^\top \mathrm{Cov}(\mathbf{\hat{\beta}}) g(z) = \sigma^2 g(z)^\top M^{-1}_{Lap}(\epsilon) g(z)
\end{eqnarray}

As it can be seen, under the LapRLS model one can obtain a D-optimal design by maximizing the determinant of $M_{Lap}(\epsilon)$ and a G-optimal design by minimizing $\operatornamewithlimits{max}\limits_{z \in \mathcal{X}} d(z,\epsilon)$. Similarly to the optimal design of experiments in Euclidean space, we prove next an equivalence theorem on Riemannian manifolds that shows how the D and G optimality criteria lead to the same optimal design. Before the equivalence theorem is discussed, we need to prove some auxiliary results. The proofs of these propositions are provided in the Appendix.

\begin{Prop}\label{prop1}
Let $\epsilon_1$ and $\epsilon_2$ be two designs with the corresponding information matrices $M_{Lap}(\epsilon_1)$ and $M_{Lap}(\epsilon_2)$. Then 
\begin{equation} 
 M_{Lap}(\epsilon_{3})=(1-\alpha) M_{Lap} (\epsilon_1) + \alpha M_{Lap} (\epsilon_2),
\end{equation}
where $M_{Lap}(\epsilon_3)$ is the information matrix of the design 
\begin{equation} 
\epsilon_{3}=(1-\alpha)\epsilon_{1}+\alpha \epsilon_2, \; \mathrm{for} \;  0<\alpha<1.
\end{equation}
\end{Prop}


\begin{Prop}\label{prop2}
Let $\epsilon_1$ and $\epsilon_2$ be two designs with the corresponding information matrices $M_{Lap}(\epsilon_1)$ and $M_{Lap}(\epsilon_2)$. Then 
\begin{equation} 
\frac{d \log |M_{Lap}(\epsilon_3)|}{d \alpha}=\Tr \Big\{ M_{Lap}^{-1}(\epsilon_3)\big [M_{Lap}(\epsilon_2)-M_{Lap}(\epsilon_1)\big] \Big\},
\end{equation}
where $M_{Lap}(\epsilon_3)$ is the information matrix of the design 
\begin{equation} 
\epsilon_3=(1-\alpha)\epsilon_1+\alpha\epsilon_2,  \; \mathrm{for} \;  0<\alpha<1.
\end{equation}
\end{Prop}


\begin{Prop} \label{prop3}
For any continuous design $\epsilon$,
\begin{enumerate}
\item \begin{equation} \label{dint}
\int_{z\in \mathcal{X}} d (z,\epsilon) \xi(z) dz=p - \Tr \Big \{ M^{-1}_{Lap}(\epsilon)  C \Big \} 
\end{equation}
\item \begin{equation} \label{maxd_ineq}
\max_{z \in \mathcal{X}} d(z,\epsilon) \geq p - \Tr \Big \{ M^{-1}_{Lap}(\epsilon) C \Big \} 
\end{equation}
\end{enumerate}
\end{Prop}


\begin{Prop}\label{prop4}
The function $\log |M_{Lap}(\epsilon)|$ is a strictly concave function.
\end{Prop}


Based on Propositions \ref{prop1}-\ref{prop4},  we can now prove the equivalence theorem for the LapRLS model. In summary, the following theorem demonstrates that the D-optimal design and G-optimal design are equivalent on the Riemannian manifold $\mathcal{M}$. It also provides the theoretical value of maximum prediction variance of the LapRLS model when the D/G optimal design is achieved. 
\begin{theorem}[\textbf{Equivalence Theorem on Manifolds}] \label{Thm1}
The following statements are equivalent:
\begin{enumerate}
\item the design $\epsilon^*$ maximizes $\det(M_{Lap}(\epsilon))$
\item the design $\epsilon^*$ minimizes $\operatornamewithlimits{max}\limits_{z \in \mathcal{X}} d(z,\epsilon)$
\item $\operatornamewithlimits{max}\limits_{z \in \mathcal{X}} d(z,\epsilon^*)=p-\Tr \Big \{ M^{-1}_{Lap}(\epsilon^*) C \Big \}$ 
\end{enumerate}
\end{theorem}

\noindent \textbf{Proof} 

\textbf{(1)} \textbf{1 $\Rightarrow$ 2}

Let $\epsilon^*$ be the design that maximizes $| M_{Lap}(\epsilon)|$ and define $\tilde{\epsilon}=(1-\alpha)\epsilon^*+\alpha \epsilon$, where $\epsilon$ is some arbitrary design.
According to Proposition \ref{prop2}, we have that
\begin{eqnarray}
\frac{d \log |M_{Lap}(\tilde{\epsilon})|}{d \alpha} &=& \Tr \Big\{ M_{Lap}^{-1}(\tilde{\epsilon})\big [M_{Lap}(\epsilon)-M_{Lap}(\epsilon^*)\big] \Big\}
\end{eqnarray}
When $\alpha=0$, we have $\tilde{\epsilon}=\epsilon^*$. Thus
\begin{eqnarray}
\frac{d \log |M_{Lap}(\tilde{\epsilon})|}{d \alpha} \bigg |_{\alpha=0} 
&=& \Tr \Big\{  M_{Lap}^{-1}(\epsilon^*) M_{Lap}(\epsilon) \Big\}- p
\end{eqnarray}
Since $\epsilon^*$ is the maximal solution, then 
\begin{eqnarray}
\Tr \Big\{  M_{Lap}^{-1}(\epsilon^*) M_{Lap}(\epsilon) \Big\}- p \leq 0.
\end{eqnarray}
Without loss of generality, assume the design $\epsilon$ has only one instance $z \in \mathcal{X}$.  Then we have
\begin{eqnarray}
M_{Lap}(\epsilon)&=&g(z)g(z)^\top+C 
\end{eqnarray}
and 
\begin{eqnarray*}
\Tr \Big\{  M_{Lap}^{-1}(\epsilon^*) \big [g(z)g(z)^\top+C \big] \Big\}- p &=& \Tr \Big\{  M_{Lap}^{-1}(\epsilon^*) g(z)g(z)^\top \Big\}+ \Tr \Big\{  M_{Lap}^{-1}(\epsilon^*) C  \Big\} - p\\
&=& d (z,\epsilon^*) + \Tr \Big\{  M_{Lap}^{-1}(\epsilon^*) C  \Big\} - p \\
&\leq& 0  \label{Mderivative}
\end{eqnarray*}
Thus 
\begin{equation} \label{dleq}
d (z,\epsilon^*) \leq p - \Tr \Big\{  M_{Lap}^{-1}(\epsilon^*) C  \Big\} 
\end{equation}
In addition, based on Proposition \ref{prop3}, we have 
\begin{equation} \label{dgeq}
\max_{z \in \mathcal{X}} d(z,\epsilon^*) \geq p - \Tr \Big \{ M^{-1}_{Lap}(\epsilon^*) C \Big \} 
\end{equation}
Combining (\ref{dleq}) and (\ref{dgeq}), we can conclude that the D-optimal design $\epsilon^*$ minimizes $\operatornamewithlimits{max}\limits_{z \in \mathcal{X}} d(z,\epsilon)$.

\rightline{$\blacksquare$}

\textbf{(2)} \textbf{2 $\Rightarrow$ 1}

Let $\epsilon^*$ be the design that minimizes $\operatornamewithlimits{max}\limits_{z \in \mathcal{X}} d(z,\epsilon)$, but assume it is not D-optimal. 
Based on Proposition \ref{prop4}, we know there must exist a design $\epsilon$ such that:
\begin{eqnarray} \label{increaseD}
\frac{d \log |(1-\alpha) M_{Lap}(\epsilon^*) + \alpha M_{Lap}(\epsilon)|}{d \alpha} \bigg |_{\alpha=0} &=& \Tr \Big\{  M_{Lap}^{-1}(\epsilon^*) M_{Lap}(\epsilon) \Big\}- p >0
\end{eqnarray}
where
\begin{eqnarray}
M_{Lap}(\epsilon) = \int_{z \in \mathcal{X}} \xi(z) g(z)g(z)^\top dz +C.
\end{eqnarray}
Then
\begin{eqnarray*}
\Tr \Big\{  M_{Lap}^{-1}(\epsilon^*) M_{Lap}(\epsilon) \Big\}- p &=& \Tr \Big\{  M_{Lap}^{-1}(\epsilon^*) \big[\int_{z \in \mathcal{X}} \xi(z) g(z)g(z)^\top dz +C \big]  \Big\}- p \\
&=& \int_{z \in \mathcal{X}} \xi(z)  \Tr \Big\{ g(z)^\top  M_{Lap}^{-1}(\epsilon^*)  g(z) \Big\} dz + \Tr \Big\{ M_{Lap}^{-1}(\epsilon^*) C \Big\}- p \\
&=& \int_{z \in \mathcal{X}} \xi(z)  d(z,\epsilon^*) dz + \Tr \Big\{ M_{Lap}^{-1}(\epsilon^*) C \Big\}- p 
\end{eqnarray*}
Since $\epsilon^*$ is the design that minimizes $\operatornamewithlimits{max}\limits_{z \in \mathcal{X}}d(z,\epsilon)$, by Proposition \ref{prop3},  we have 
\begin{eqnarray}
\max_{z \in \mathcal{X}} d(z,\epsilon^*)= p - \Tr \Big \{ M^{-1}_{Lap}(\epsilon^*) C \Big \} 
\end{eqnarray}
Thus, for any $z \in \mathcal{X}$, 
\begin{eqnarray}
 d(z,\epsilon^*)  &\leq& p - \Tr \Big\{  M_{Lap}^{-1}(\epsilon^*) C \Big\} 
\end{eqnarray}
Then
\begin{eqnarray}
\int_{z \in \mathcal{X}} \xi(z)  d(z,\epsilon^*) dz &\leq& \int_{z \in \mathcal{X}}  \xi(z) \Bigg(p-\Tr \Big\{ M_{Lap}^{-1}(\epsilon^*) C \Big\} \Bigg) dz \\
&=& p-\Tr \Big\{ M_{Lap}^{-1}(\epsilon^*) C \Big\}
\end{eqnarray}
Therefore, we have
\begin{eqnarray*}
\Tr \Big\{  M_{Lap}^{-1}(\epsilon^*) M_{Lap}(\epsilon) \Big\}- p &=& \int_{z \in \mathcal{X}} \xi(z)  d(z,\epsilon^*) dz + \Tr \Big\{ M_{Lap}^{-1}(\epsilon^*) C \Big\}- p \\
&\leq& p - \Tr \Big\{  M_{Lap}^{-1}(\epsilon^*) C \Big\} + \Tr \Big\{ M_{Lap}^{-1}(\epsilon^*) C \Big\}- p \\
&=& 0.
\end{eqnarray*}
This contradicts with (\ref{increaseD}). Therefore, the design $\epsilon^*$ is also D-optimal. 

\rightline{$\blacksquare$}

\textbf{(3) 1 $\Rightarrow$ 3}

Let $\epsilon^*$ be the D-optimal design. From the previous proof, in particular Equation (\ref{dleq}), we know that 
\begin{eqnarray}
\max_{z\in \mathcal{X}} d(z,\epsilon^*) = 
p - \Tr \Big \{ M^{-1}_{Lap}(\epsilon^*) C \Big \}. 
\end{eqnarray}
\rightline{$\blacksquare$}

\textbf{(4) 3 $\Rightarrow$ 1}

Let $\epsilon^*$ be the design such that 
\begin{eqnarray}
\max_{z\in \mathcal{X}} d(z,\epsilon^*) = p - \Tr \Big \{ M^{-1}_{Lap}(\epsilon^*) C \Big \}. 
\end{eqnarray}
Then, for any $z \in \mathcal{X}$, 
\begin{equation} \label{eq65}
d (z,\epsilon^*) + \Tr \Big\{  M_{Lap}^{-1}(\epsilon^*) C  \Big\} - p \leq 0.
\end{equation}
Based on the previous proof, we know that equation (\ref{eq65}) implies that there is no improving direction for the D-optimal criteria. Thus $\epsilon^*$ is the D-optimal design. \\
\rightline{$\blacksquare$}

\textbf{(5)} Since 1 $\Leftrightarrow$ 2, 1 $\Leftrightarrow$ 3, then \textbf{2 $\Leftrightarrow$ 3} and the equivalence theoreme is proved. 

\rightline{$\blacksquare$} 

Different from the classical equivalence theorem on Euclidean space, 
Theorem \ref{Thm1} demonstrates the equivalence of D-optimal design and G-optimal design on the Riemannian manifold. In addition, for any given design $\epsilon$, Equation (\ref{maxd_ineq}) provides a new lower bound for the maximum prediction variance.  Theorem \ref{Thm1} shows that this lower bound (\ref{maxd_ineq}) can be achieved at the D/G optimal design $\epsilon^*$. Therefore, Theorem \ref{Thm1} also provides a theoretical justification that the optimal D/G design $\epsilon^*$  minimizes the maximum prediction variance of the model.


\section{Proposed Algorithm and Convergence Analysis}\label{sec:alg}
Before we discuss the proposed algorithm for finding optimal experimental design on manifolds, some auxiliary results need to be given, whose proofs are in the Appendix.  

\vspace{0.1cm}
\begin{Prop}\label{prop5}
Let $M_{Lap}(\epsilon_k)$ be the information matrix of the design $\epsilon_k$ at $k$-th iteration. Let $M_{Lap}(\epsilon(z))$ be the information matrix of the design concentrated at one single point $z$. Given $\epsilon_{k+1}=(1-\alpha)\epsilon_{k}+\alpha \epsilon(z)$, then  
\begin{eqnarray} \label{eqprop4}
|M_{Lap}(\epsilon_{k+1})|=(1-\alpha)^p \Big|M_{Lap}(\epsilon_k) \Big| \Big[1+\frac{\alpha}{1-\alpha} d(z,\epsilon_k)+\frac{\alpha}{1-\alpha} \Tr(M^{-1}_{Lap}(\epsilon_k) C) \Big] 
\end{eqnarray}
\end{Prop}

\begin{Prop}\label{prop6}
Let $M_{Lap}(\epsilon_k)$ be the information matrix of the design $\epsilon_k$ at $k$-th iteration. Construct the design $\epsilon_{k+1}$ at $(k+1)$-th iteration as 
\begin{eqnarray} 
\epsilon_{k+1}=(1-\alpha_k)\epsilon_k+\alpha_k \epsilon(z_{k+1})
\end{eqnarray}
where 
\begin{eqnarray} 
0 < \alpha_k \leq \frac{d(z_{k+1},\epsilon_k)-(p-\Tr(M^{-1}_{Lap}(\epsilon_k) C))}{p[d(z_{k+1},\epsilon_k)-(1-\Tr(M^{-1}_{Lap}(\epsilon_k) C))]}, \;  \; z_{k+1}=\operatornamewithlimits{argmax}\limits_{z \in \mathcal{X}} d (z, \epsilon_k).
\end{eqnarray}
Then  the resulting sequence $\Big\{|M_{Lap}(\epsilon_k)|\Big\}_k$ is monotonic increasing. 
\end{Prop}

Based on Propositions \ref{prop5} and \ref{prop6}, the new algorithm for finding a D-G optimal experimental design on a manifold is shown in Algorithm 1. Note how after obtaining an optimal design for the data to be labeled, and obtaining the corresponding labels, we can use both labeled {\em and} unlabeled instances to train the manifold regularized model (\ref{eqn:mrobj2}). 

\begin{algorithm}[!htbp]
\caption{Optimal Design of Experiments on Manifolds ( {ODOEM})}
\begin{algorithmic}
\Statex \textbf{Input:} Some initial design $\epsilon_k$,
\begin{eqnarray*}
\epsilon_k=\left \{   \begin{array}{cccc}
                  z_1, z_2, ..., z_k\\
                  q_1,q_2,...,q_k
                \end{array}
              \right\}, \; \mathrm{where} \; \sum_{i=1}^k q_i=1
\end{eqnarray*}

Compute the information matrix
\begin{eqnarray}
M_{Lap}(\epsilon_k)=\sum_{i=1}^k q_i g(z_i)g(z_i)^\top +C
\end{eqnarray}

\While{optimal design is not achieved}
\begin{enumerate}
\item Find $z_{k+1}$ s.t.
\begin{eqnarray}
z_{k+1}=\operatornamewithlimits{argmax}\limits_{z \in \mathcal{X}} d(z,\epsilon_k) 
\end{eqnarray}
\item Update the design 
\begin{eqnarray}
\epsilon_{k+1}=(1-\alpha_k)\epsilon_k+\alpha_k \epsilon(z_{k+1})
\end{eqnarray}
where $\alpha_k$ is a user choice that satisfies 
\begin{eqnarray}
0 < \alpha_k \leq \frac{d(z_{k+1},\epsilon_k)-[p-\Tr(M^{-1}_{Lap}(\epsilon_k) C)]}{p\{d(z_{k+1},\epsilon_k)-[1-\Tr(M^{-1}_{Lap}(\epsilon_k) C)]\}} 
\end{eqnarray}
\item Compute the information matrix $M_{Lap}(\epsilon_{k+1})$, set $k=k+1$ and repeat step 1-3.
\end{enumerate}
\EndWhile
\Statex \textbf{Output:} Optimal design on manifolds $\epsilon^*$.
\end{algorithmic}
\end{algorithm}

We next provide a convergence analysis of the proposed algorithm.

\begin{theorem}[\textbf{Convergence Theorem}] 
The iterative procedure in Algorithm 1 converges to the D-optimal design $\epsilon^*$,
\begin{eqnarray}
\lim_{k \rightarrow \infty} |M_{Lap}(\epsilon_k)|=|M_{Lap}(\epsilon^*)| 
\end{eqnarray}
\end{theorem}

\noindent \textbf{Proof} 

Let the design $\epsilon_0$ not be D-optimal. Based on Proposition \ref{prop6}, we have
\begin{eqnarray}
|M_{Lap}(\epsilon_0)| < |M_{Lap}(\epsilon_1)| < \cdots < |M_{Lap}(\epsilon_k)| < \cdots \leq |M_{Lap}(\epsilon^*)| 
\end{eqnarray}
It is known that any bounded monotone sequence converges. Thus the sequence $|M_{Lap}(\epsilon_0)|$, $|M_{Lap}(\epsilon_1)|$, ..., $|M_{Lap}(\epsilon_k)|$ converges to some limit $|M_{Lap}(\hat\epsilon)|$. Next we need to show 
\begin{eqnarray}
|M_{Lap}(\hat\epsilon)|=|M_{Lap}(\epsilon^*)| 
\end{eqnarray}
The proof proceeds by contradiction. Assume
\begin{eqnarray} \label{eqassumption}
|M_{Lap}(\hat\epsilon)|<|M_{Lap}(\epsilon^*)| 
\end{eqnarray}
By the convergence of the sequence $|M_{Lap}(\epsilon_0)|$, $|M_{Lap}(\epsilon_1)|$, ..., $|M_{Lap}(\epsilon_k)|$, we know that, for $\forall \eta >0$, there $\exists k_0 \in \mathbb{N}$ s.t. 
\begin{eqnarray}
|M_{Lap}(\epsilon_{k+1})|-|M_{Lap}(\epsilon_k)| < \eta \; \mathrm{for} \; \forall k > k_0 
\end{eqnarray}
Based on Proposition \ref{prop5}, we have
\begin{eqnarray*} 
(1-\alpha_k)^p\Big(1+\frac{\alpha_k}{1-\alpha_k} d(z_{k+1},\epsilon_k)+\frac{\alpha_k}{1-\alpha_k} \Tr(M^{-1}_{Lap}(\epsilon_k) C) \Big)  |M_{Lap}(\epsilon_k)| -|M_{Lap}(\epsilon_k)| &<& \eta 
\end{eqnarray*}
Then,
\begin{eqnarray}\label{eq93} 
(1-\alpha_k)^p\Big(1+\frac{\alpha_k}{1-\alpha_k} [d(z_{k+1},\epsilon_k)+\Tr(M^{-1}_{Lap}(\epsilon_k) C) ] \Big) & <& 1+ \eta |M_{Lap}(\epsilon_k)|^{-1}
\end{eqnarray}
Defining $\tau_k=d(z_{k+1},\epsilon_k)-[p-\Tr(M^{-1}_{Lap}(\epsilon_k) C)]$,  we can rewrite (\ref{eq93}) as 
\begin{eqnarray} \label{eq94}
(1-\alpha_k)^p\Big(1+\frac{\alpha_k}{1-\alpha_k} [\tau_k+p] \Big)  < 1+ \eta |M_{Lap}(\epsilon_k)|^{-1}
\end{eqnarray}
Next, define a function $\mathrm{T}(\tau_k, \alpha_k)$ as
\begin{eqnarray}\label{eq:tau:def}
\mathrm{T}(\tau_k, \alpha_k)=(1-\alpha_k)^p\Big(1+\frac{\alpha_k}{1-\alpha_k} [\tau_k+p] \Big) 
\end{eqnarray}
such that
\begin{eqnarray}
\frac{\partial \mathrm{T}}{\partial \tau_k}&=&(1-\alpha_k)^p \frac{\alpha_k}{1-\alpha_k} 
\end{eqnarray}
Clearly, $\frac{\partial \mathrm{T}}{\partial\tau_k} > 0$ for $0 < \alpha_k<1$. Thus, for a given $0 < \alpha_k<1$, $\mathrm{T}(\tau_k, \alpha_k)$  is a monotonic increasing function with respect to $\tau_k$. On the other hand, 
\begin{eqnarray}
\frac{\partial \mathrm{T}}{\partial\alpha_k}=-p(1-\alpha_k)^{p-1}\Big(1+\frac{\alpha_k}{1-\alpha_k} [\tau_k+p] \Big) +(1-\alpha_k)^p [\tau_k+p] \frac{1}{(1-\alpha_k)^2} 
\end{eqnarray}
Let $\frac{\partial \mathrm{T}}{\partial\alpha_k} \geq 0$, we have
\begin{eqnarray}
(1-\alpha_k)^{p-2} [\tau_k+p]  &\geq& p(1-\alpha_k)^{p-1}\Big(1+\frac{\alpha_k}{1-\alpha_k} [\tau_k+p] \Big) \\
 \alpha_k &\leq&\frac{\tau_k}{p(p + \tau_k-1)}
\end{eqnarray}
Thus, for $0 < \alpha_k \leq \frac{\tau_k}{p(p + \tau_k-1)}$ and  $\tau_k>0$, $\mathrm{T}(\tau_k,\alpha_k)$ is a monotone increasing function. In particular, plugging in the expression for $\tau_k$, we get
\begin{eqnarray}
\frac{\tau_k}{p(p + \tau_k-1)}=\frac{d(z_{k+1},\epsilon_k)-[p-\Tr(M^{-1}_{Lap}(\epsilon_k) C)]}{p(d(z_{k+1},\epsilon_k)-[1-\Tr(M^{-1}_{Lap}(\epsilon_k) C)])}. 
\end{eqnarray}
Notice that $0 < \alpha_k \leq \frac{\tau_k}{p(p + \tau_k-1)}$ is the same choice of $\alpha_k$ in the proposed Algorithm 1. 

From the assumption (\ref{eqassumption}), Proposition \ref{prop3} and Theorem \ref{Thm1}, it follows that $\tau_k>0$. This guarantees the existence of $\alpha_k$ such that $0 < \alpha_k \leq \frac{\tau_k}{p(p + \tau_k-1)}$. Thus, for any $\tau_k>0$ and $0 < \alpha_k \leq \frac{\tau_k}{p(p + \tau_k-1)}$, we have $\mathrm{T}(\tau_k,\alpha_k) > 1$. Note that $\eta$ is an arbitrary positive number in equation (\ref{eq94}), which implies $\tau_k$ need to be an infinitely small positive number to satisfy equation (\ref{eq94}), i.e. given $\forall\zeta>0$, there $\exists \tilde{k} (\zeta) \in \mathbb{N}$ s.t.
\begin{eqnarray}
\tau_k=d(z_{k+1},\epsilon_k)-[p-\Tr(M^{-1}_{Lap}(\epsilon_k) C)] < \zeta \; \quad \mathrm{for} \; k>\tilde{k}(\zeta) 
\end{eqnarray}
However, based on the assumption (\ref{eqassumption}) and Theorem \ref{Thm1}, we have that
\begin{eqnarray}
d(z_{k+1},\epsilon_k)-[p-\Tr(M^{-1}_{Lap}(\epsilon_k) C)] \geq \delta_k >0 \; \quad \mathrm{for} \; \forall k.
\end{eqnarray}
Choosing $\zeta < \delta_k$, we have a contradiction, and therefore, the convergence theorem is proved.

\rightline{$\blacksquare$}

From the derivation of Algorithm 1, it is not difficult to notice that ODOEM is a model-dependent design. The corresponding manifold regularization model (\ref{eqn:LapRLS}) need to be trained after a desired number of instances is labeled. As it is shown before, Algorithm 1 is a converging algorithm on a continuous design space. However, sometimes the experimental design space is not continuous and only a set of candidate points is available. For a discrete design space with a set of candidate points, and in analogy to ODOE's on Euclidean spaces, one can evaluate each candidate point and choose the point with maximum prediction variance. The resulting sequence of $|M_{Lap}(\epsilon_k)|$ is still monotonic increasing, since
\begin{eqnarray*}
|M_{Lap}(\epsilon_k)+g(z_{k+1})g(z_{k+1})^\top|&=&|M_{Lap}(\epsilon_k)|[1+g(z_{k+1})^\top M^{-1}_{Lap}(\epsilon_k) g(z_{k+1})] \\
&>& |M_{Lap}(\epsilon_k)|
\end{eqnarray*}
where $z_{k+1}=\operatornamewithlimits{argmax}\limits_{z \in \mathcal{X}\setminus Z_k} d(z,\epsilon_k)=\operatornamewithlimits{argmax}\limits_{z \in \mathcal{X}\setminus Z_k} g(z)^\top M^{-1}_{Lap}(\epsilon_k) g(z)$.


\section{Numerical Results}\label{sec:numerical}

To illustrate the empirical performance of the proposed ODOEM algorithm in practice, we consider its application to both synthetic datasets  and also its application to the high dimensional real-world image datasets. The synthetic datasets are low dimensional manifold examples that permit straightforward visualization of the resulting designs and are shown first.

\subsection{Synthetic Manifold Datasets}
\begin{figure}[!htbp]
        \centering
        \begin{subfigure}[b]{0.31\textwidth}
            \centering
            \includegraphics[width=\textwidth]{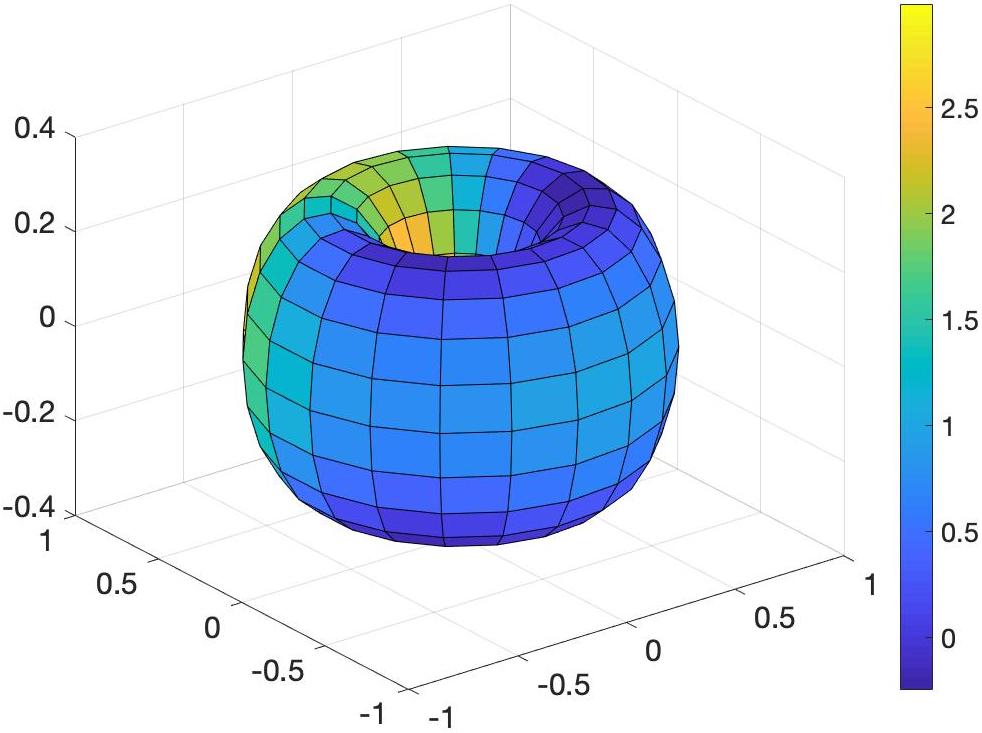}
        \end{subfigure}
        \quad
        \begin{subfigure}[b]{0.31\textwidth}
            \centering
            \includegraphics[width=\textwidth]{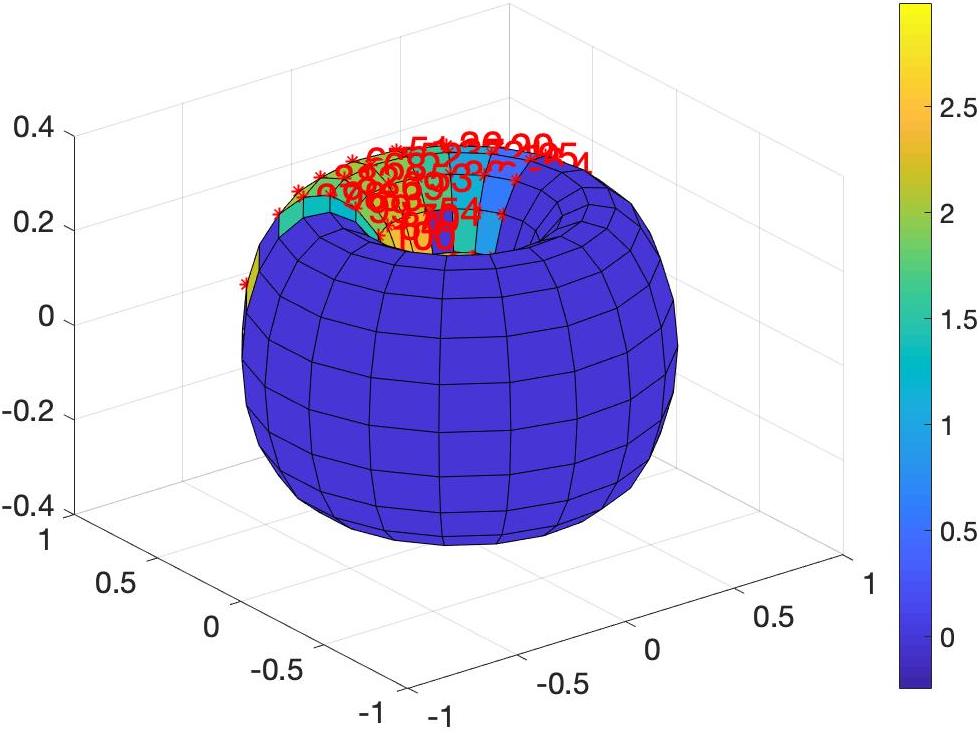}
        \end{subfigure}
        \quad
        \begin{subfigure}[b]{0.31\textwidth}
            \centering
            \includegraphics[width=\textwidth]{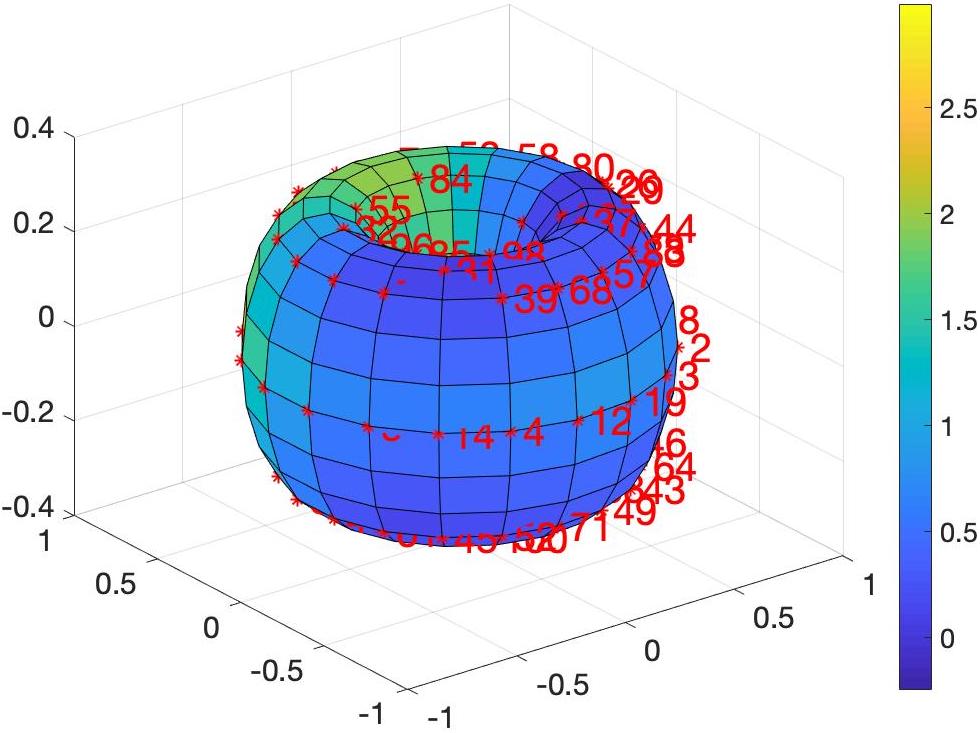}
        \end{subfigure}
        \vskip\baselineskip
        \begin{subfigure}[b]{0.31\textwidth}
            \centering
            \includegraphics[width=\textwidth]{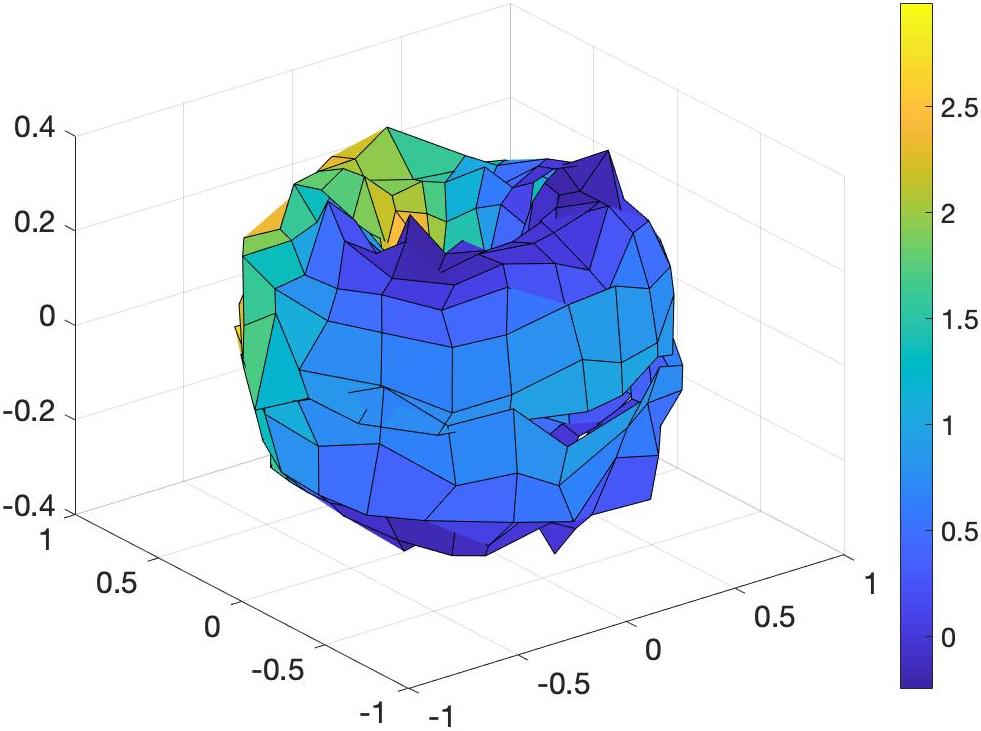}
            \caption{True Function}
        \end{subfigure}
        \quad
        \begin{subfigure}[b]{0.31\textwidth}
            \centering
            \includegraphics[width=\textwidth]{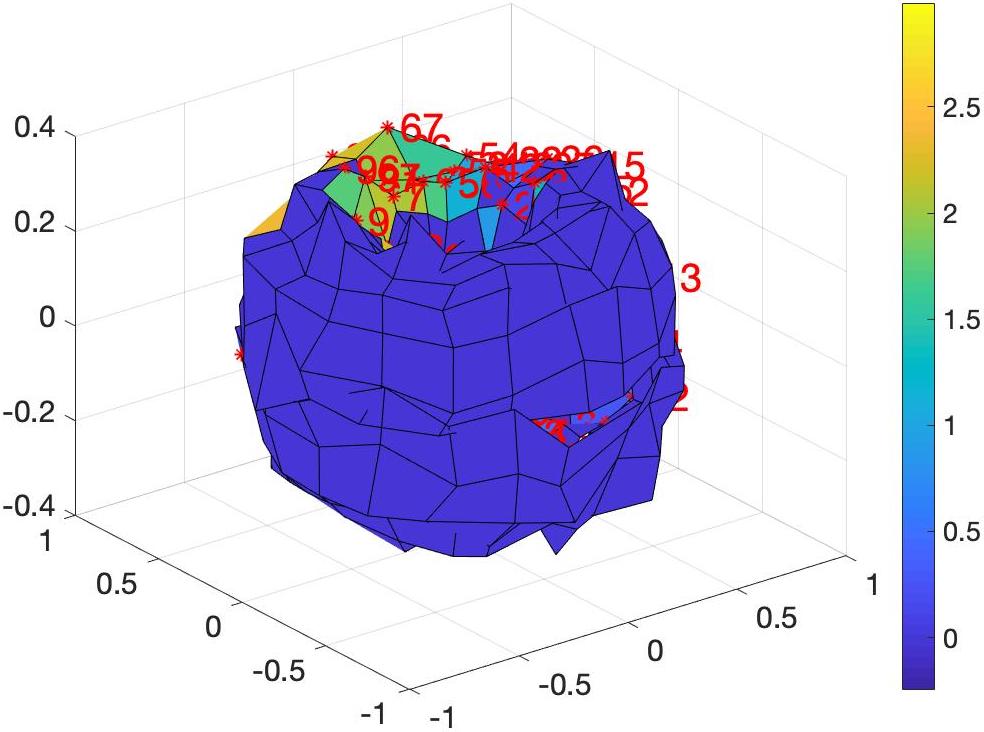}
            \caption{Classical D-optimal}
        \end{subfigure}
        \quad
        \begin{subfigure}[b]{0.31\textwidth}
            \centering
            \includegraphics[width=\textwidth]{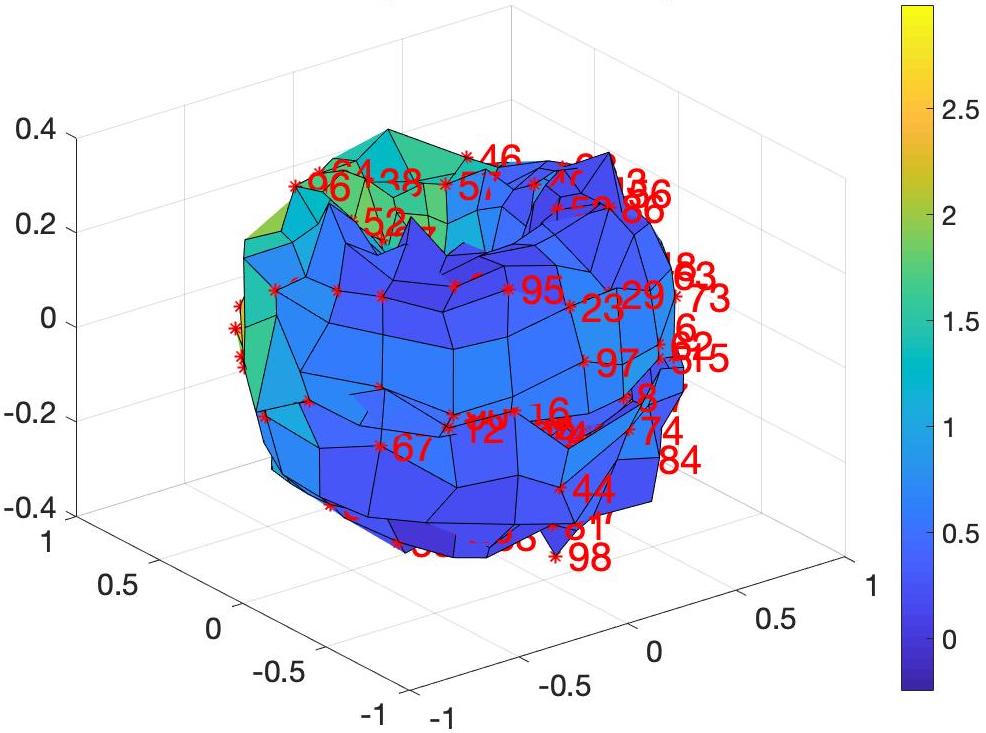}
            \caption{ODOEM}
        \end{subfigure}
        \caption{Torus example. Top: when $\{\bx_i\}_{i=1}^n$  lie on a Torus. Bottom: when $\{\bx_i\}_{i=1}^n$ are not exactly on a Torus due to noise. The simulated isotropic noise follows a normal distribution with zero mean  and variance equal to 0.03 in each ambient dimension.   (a) The colors represent the true response values defined on the Torus. (b) 100 labeled instances (red numbers) and fitted response values (colors on the surface of the Torus) by a kernel regression with the D-optimal Design. (c) 100 labeled instances (red numbers) and fitted response values  (colors on the surface of the Torus)  by ODOEM. As it can be seen, the fitted function in (c) approximates the true function on the Torus in (a) better than the fitted function in (b), with or without noise.}
        \label{fig:torus}
\end{figure}

\begin{figure}[H]
        \centering
        \begin{subfigure}[b]{0.31\textwidth}
            \centering
            \includegraphics[width=\textwidth]{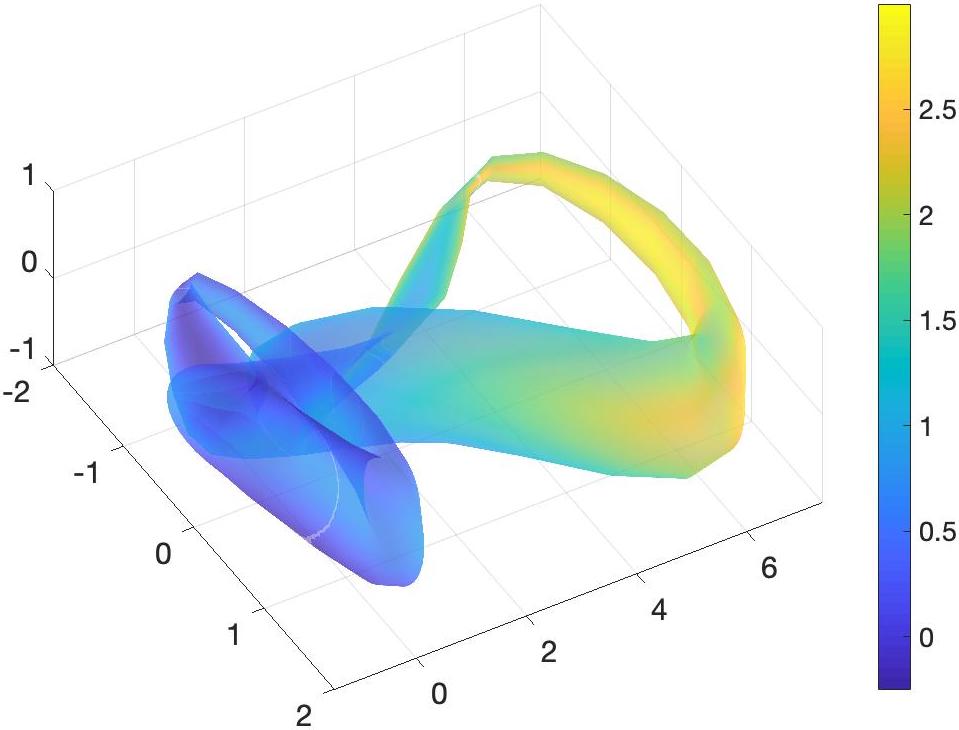}
        \end{subfigure}
        \quad
        \begin{subfigure}[b]{0.31\textwidth}
            \centering
            \includegraphics[width=\textwidth]{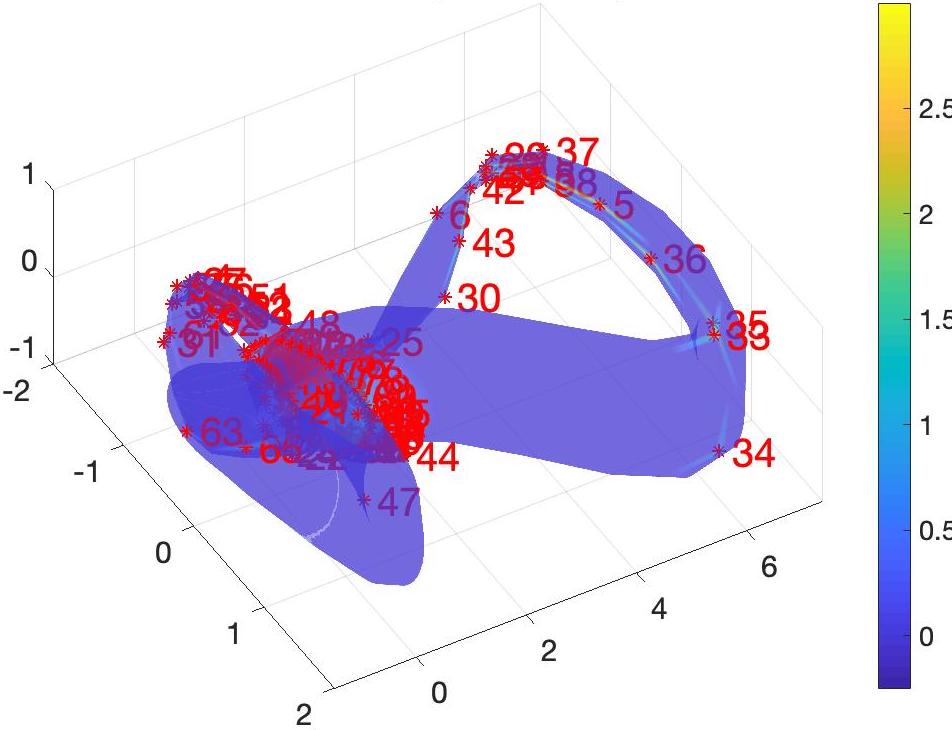}
        \end{subfigure}
        \quad
        \begin{subfigure}[b]{0.31\textwidth}
            \centering
            \includegraphics[width=\textwidth]{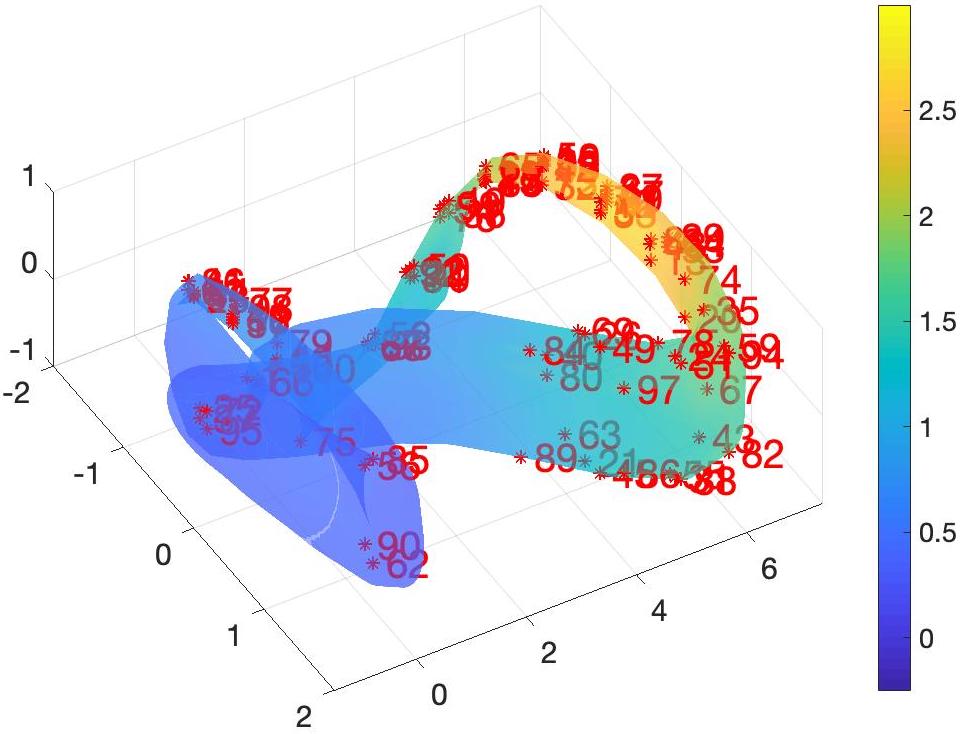}
        \end{subfigure}
         \vskip\baselineskip
         \begin{subfigure}[b]{0.31\textwidth}
            \centering
            \includegraphics[width=\textwidth]{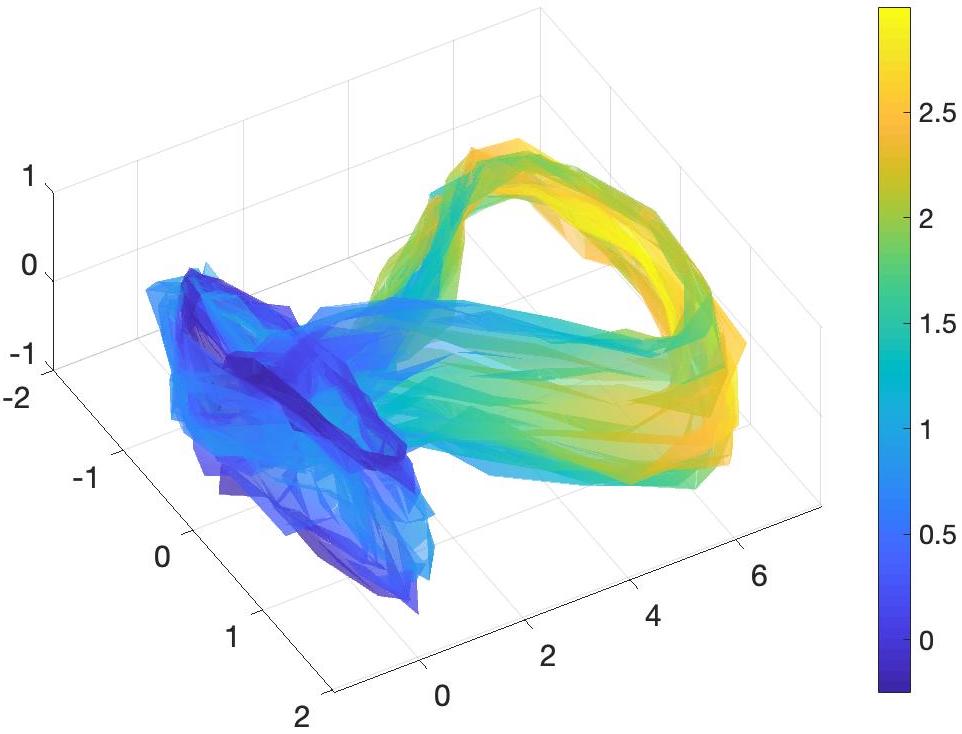}
            \caption{True Function}
        \end{subfigure}
        \quad
        \begin{subfigure}[b]{0.31\textwidth}
            \centering
            \includegraphics[width=\textwidth]{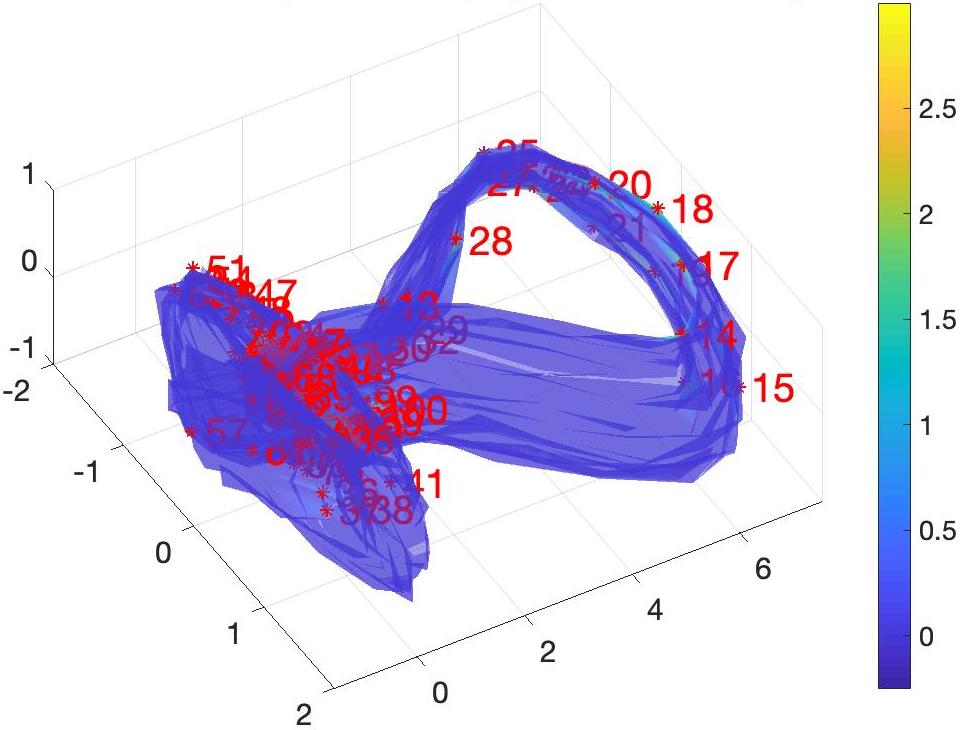}
            \caption{Classical D-optimal}
        \end{subfigure}
        \quad
        \begin{subfigure}[b]{0.31\textwidth}
            \centering
            \includegraphics[width=\textwidth]{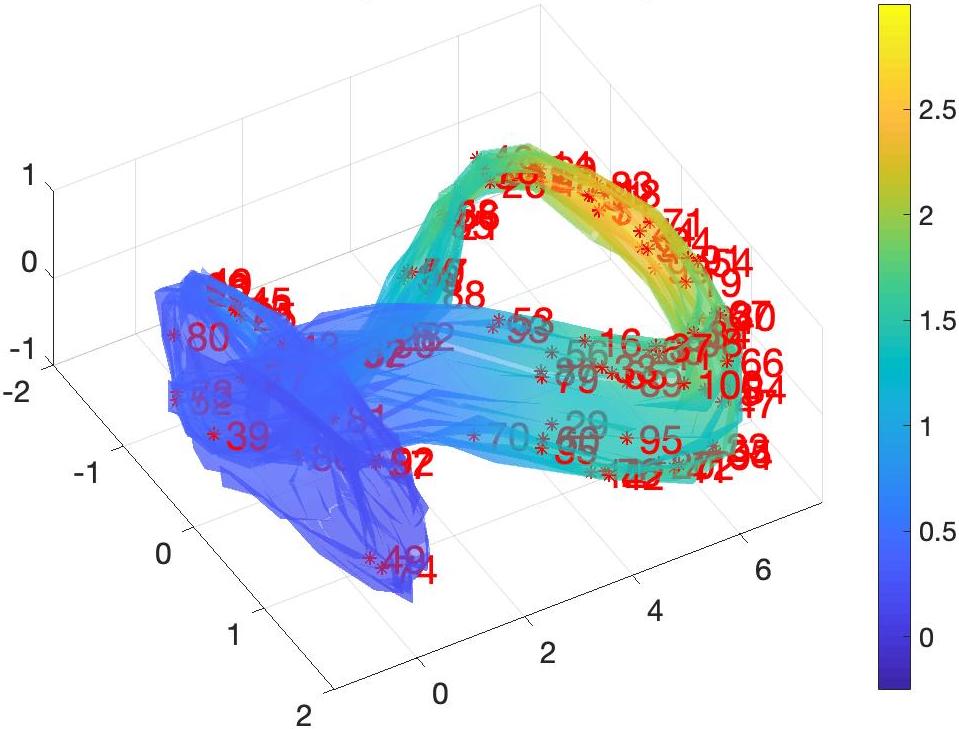}
            \caption{ODOEM}
        \end{subfigure}
        \caption{Bottle Shape of Klein Bottle example. Top: when $\{\bx_i\}_{i=1}^n$  lie on a Klein Bottle. Bottom: when $\{\bx_i\}_{i=1}^n$ are not exactly on a Klein Bottle due to noise. The simulated isotropic noise follows a normal distribution with zero mean  and variance equals to 0.06 on each coordinate on the ambient space.   (a) The colors represent the true response values defined on the Klein Bottle. (b) 100 labeled instances (red numbers) and fitted response values (colors on the surface of the Klein Bottle) by a kernel regression with the D-optimal Design. (c) 100 labeled instances (red numbers) and fitted response values  (colors on the surface of the Klein Bottle)  by ODOEM. Again, the fitted function in (c) approximates the true function on the Klein Bottle in (a) better than the fitted function in (b), with or without noise.}
        \label{fig:bottle}
\end{figure}

We generate four different two-dimensional manifold datasets: data on a Torus, on a M\"{o}bius Strip, on a figure ``8''  immersion of a Klein bottle and on a classic Klein bottle \citep{10.5555/1213748}. Each of the first three datasets contains 400 instances and the last dataset contains 1600 instances. For all four datasets, we plot these two-dimensional manifolds in a three-dimensional Euclidean space, as shown in Figures \ref{fig:torus}-\ref{fig:bottle}. The colors on these manifolds represent the corresponding response values $\{y_i\}_{i=1}^n$ or their estimates $\{\hat{y}_i\}_{i=1}^n$ based on different experimental designs. The true response values $\{y_i\}_{i=1}^n$ are defined by
\begin{equation}
    y=sin(u)+sin^2(u)+cos^2(v)
\end{equation}
where $u \in [0,2\pi)$ and $v \in [0,2\pi)$. The red numbers on the manifolds represent the sequence of labeled instances by different design algorithms. 

The regularization parameters $\lambda_A$ and $\lambda_I$ are usually selected by cross-validation. However, ODOEM is a sequential design algorithm and the order in which instances (points on the manifold) are  labeled is important. The cross-validation idea, which randomly divides the labeled instances into a training set and a validation set, is impractical in a sequential design. Thus,  we set fixed values for $\lambda_A$ and $\lambda_I$ using in our experiments $\lambda_A=0.01$ for numerical stability and generate the decreasing sequence $\lambda_I=-\ln(k/n)$, where $k$ is the number of labeled instance at the $k$-th iteration and $n$ is the total number of instances. The reason we choose a decreasing sequence of $\lambda_I$ comes from the penalized loss function (\ref{eqn:LapRLS}) and the performance evaluation criterion MSE$=\sum_{i=1}^n (y_i -\hat{f}(z_i))^2$. For the manifold regularization model, the estimated learning function $\hat{f}$ is obtained by minimizing the objective function (\ref{eqn:LapRLS}). At early iterations, there are only few labeled instances, and $\hat{f}$ would benefit more from penalizing the learning function along the manifold structure (second regularization term). As the number of labeled instances increases, larger $\lambda_I$ might not lead to smaller MSE. For example, we consider the extreme scenario when all the instances have been labeled, i.e. $k=n$. If one desires to achieve a smaller MSE$=\sum_{i=1}^n (y_i -\hat{f}(z_i))^2$, it is better to estimate $\hat{f}$ by
\begin{equation}
\hat{f}=\operatornamewithlimits{argmin}\limits_{f\in \mathcal{H}_K} \sum_{i=1}^n (y_i -f(z_i))^2,
\end{equation}
instead of using (\ref{eqn:LapRLS}). Therefore, we set $\lambda_I=-\ln(k/n)$ so that we can get a decreasing sequence of $\lambda_I$'s as $k$ increases and $\lambda_I=0$ when all the instances have been labeled. 

We compare the ODOEM algorithm with a classical D-optimal design algorithm on a kernel regression model, which does not consider the manifold structure. For both of the learning models, we choose a RBF kernel and set the range parameter to be 0.01. 

For some applications, the data may not strictly lie on a given manifold due to noise. In order to explore the robustness of the ODOEM algorithm to noise, we also let the four synthetic datasets fluctuate around their manifolds by adding noise to $\{\bx_i\}_{i=1}^n$.  In other words, for each of the four manifolds, we investigate both the case when the data $\{\bx_i\}_{i=1}^n$  lie exactly on the given manifold and the case when $\{\bx_i\}_{i=1}^n$ are not exactly on the manifold. The results are shown in Figure \ref{fig:torus}-\ref{fig:bottle}. 

Based on these results, the following comments can be made: (a) in the Torus,  M\"{o}bius Strip, and Figure “8” Immersion examples, the instances selected by the classical D-optimal design tend to be clustered in certain regions, while the instances selected by ODOEM are widely spread over these manifolds. Although the fitted values from using the classical D-optimal design is close to the true values in some very small regions on manifolds, it is clear that ODOEM provides better overall fitting performance. (b) In the Klein Bottle example, the classical D-optimal design selects relatively dispersive instances, but the function is still poorly fitted on the Bottle. It is illustrated that the kernel regression is not able to capture the manifold structure and incorporate it into the learning process. (c) ODOEM is adaptive to various manifold structures. It picks instances all over the manifolds and provides stable and superior fitting performance. In summary, on all four synthetic manifold datasets, ODOEM performs much better than kernel regression D-optimal Design in terms of instance selection and function fitting, under both of the noise-free cases and the noisy cases.

\begin{figure}[H]
        \centering
        \begin{subfigure}[b]{0.31\textwidth}
            \centering
            \includegraphics[width=\textwidth]{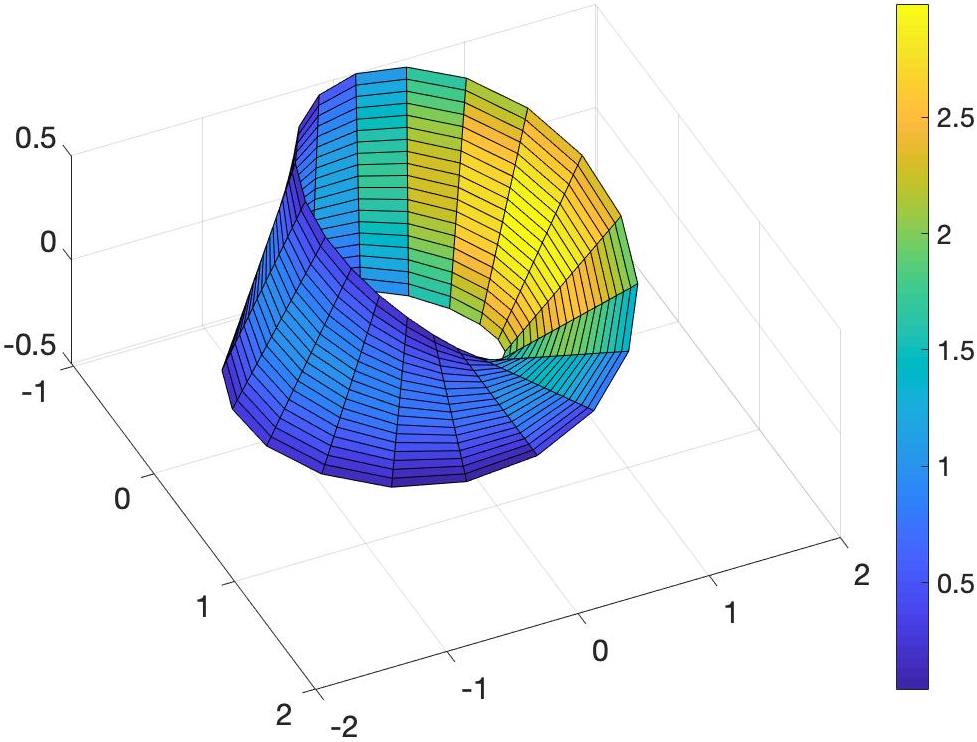}
        \end{subfigure}
        \quad
        \begin{subfigure}[b]{0.31\textwidth}
            \centering
            \includegraphics[width=\textwidth]{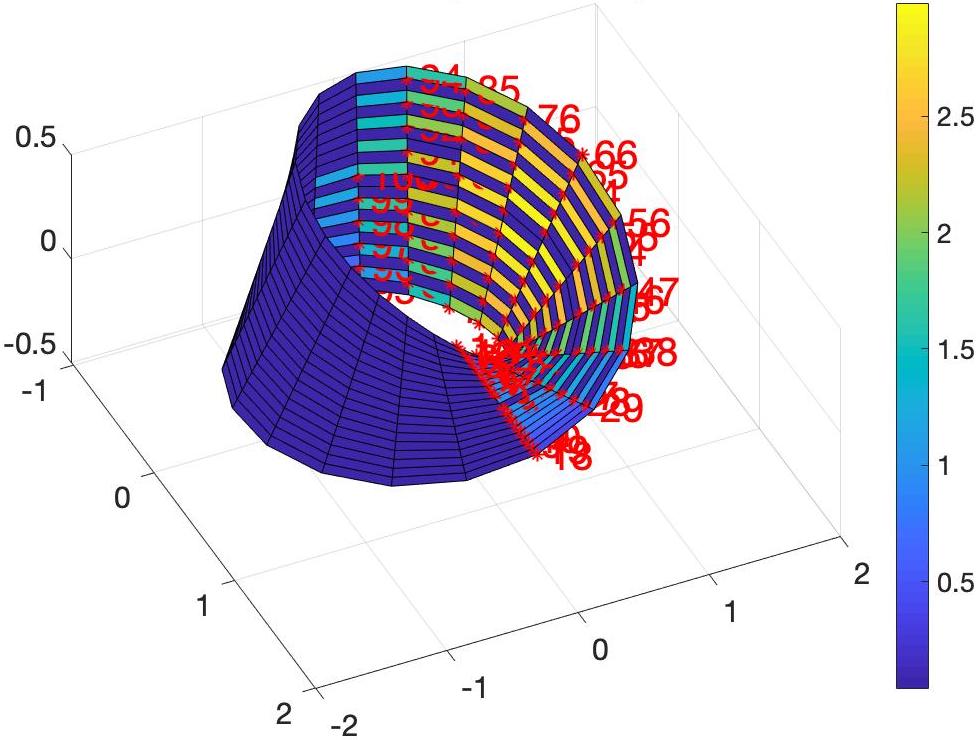}
        \end{subfigure}
        \quad
        \begin{subfigure}[b]{0.31\textwidth}
            \centering
            \includegraphics[width=\textwidth]{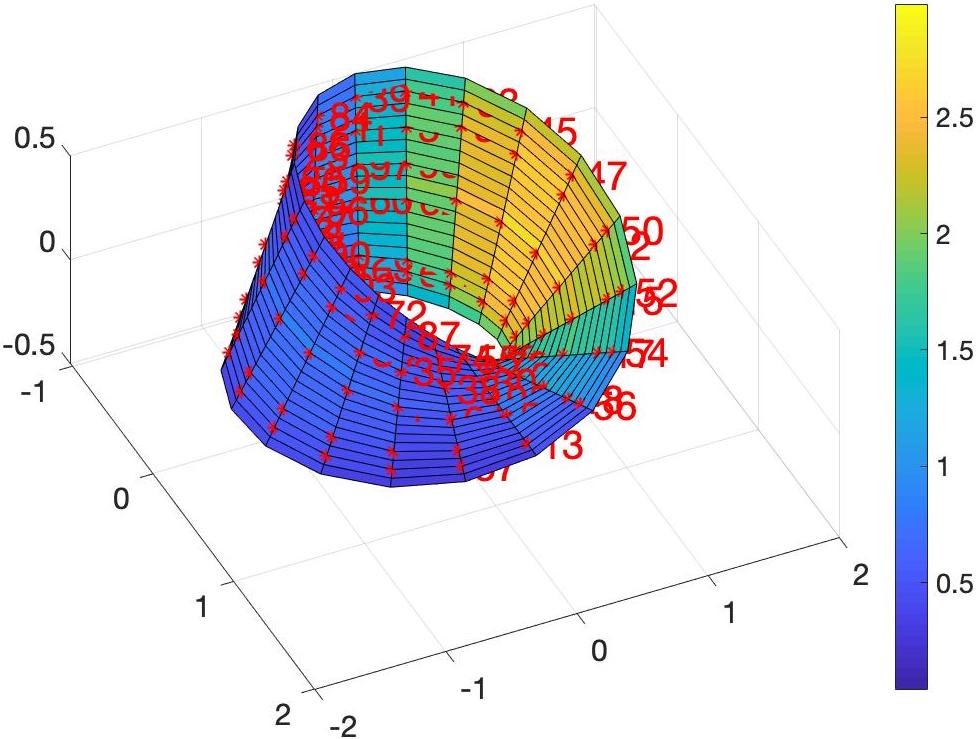}
        \end{subfigure}
         \vskip\baselineskip
        \begin{subfigure}[b]{0.31\textwidth}
            \centering
            \includegraphics[width=\textwidth]{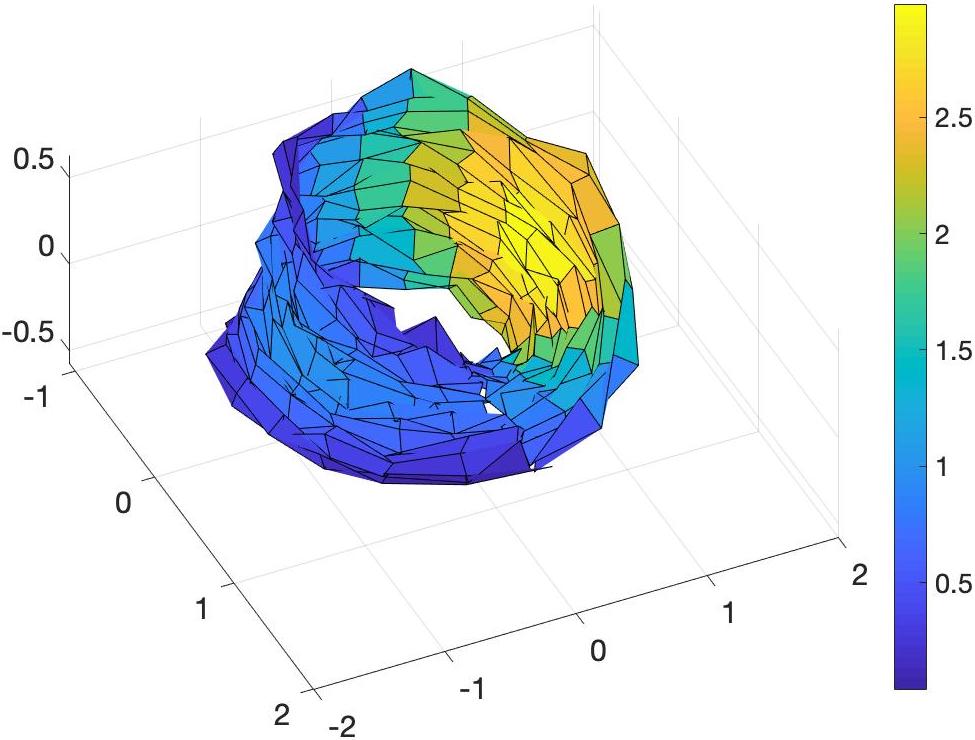}
            \caption{True Function}
        \end{subfigure}
        \quad
        \begin{subfigure}[b]{0.31\textwidth}
            \centering
            \includegraphics[width=\textwidth]{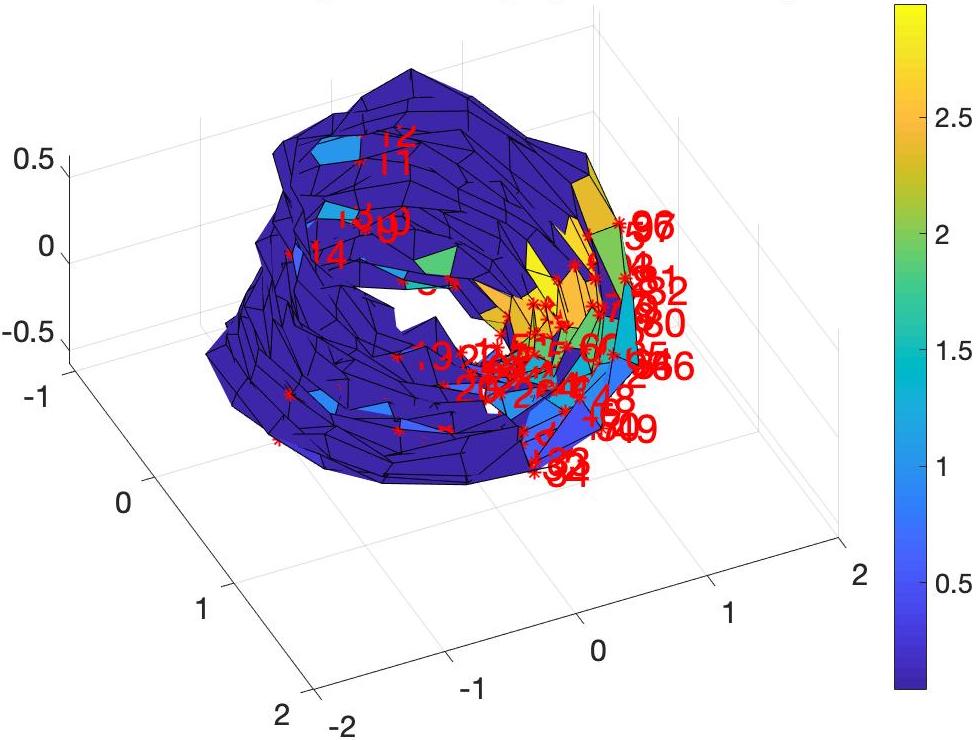}
            \caption{Classical D-optimal}
        \end{subfigure}
        \quad
        \begin{subfigure}[b]{0.31\textwidth}
            \centering
            \includegraphics[width=\textwidth]{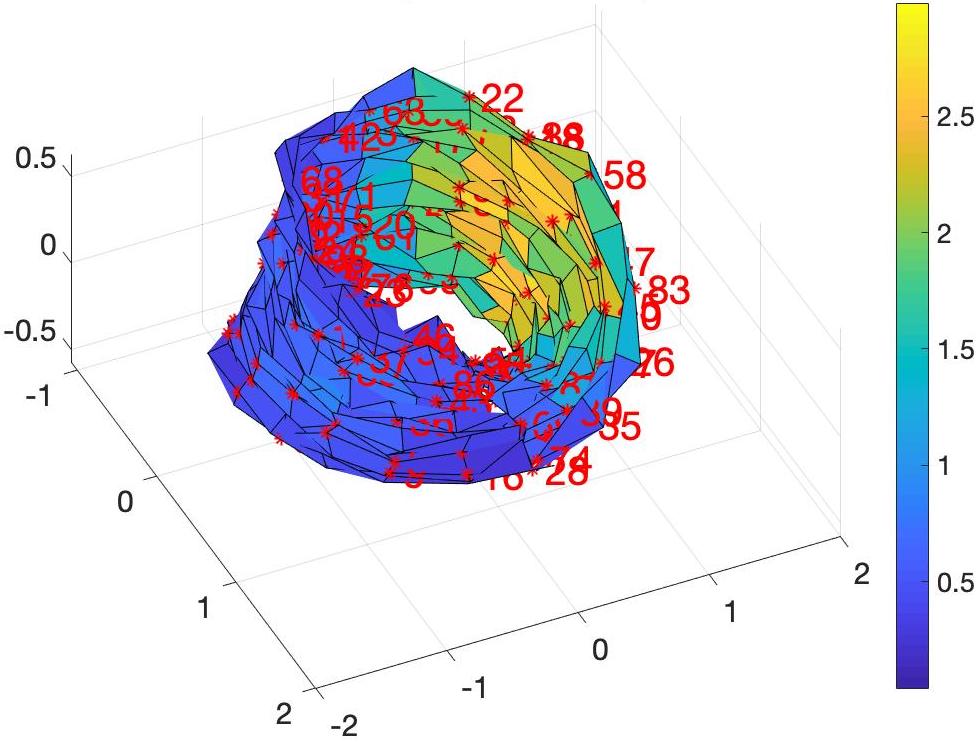}
            \caption{ODOEM}
        \end{subfigure}
         \caption{M\"{o}bius Strip example. Top: when $\{\bx_i\}_{i=1}^n$  lie on a M\"{o}bius Strip. Bottom: when $\{\bx_i\}_{i=1}^n$ are not exactly on a M\"{o}bius Strip due to noise. The simulated isotropic noise follows a normal distribution with zero mean  and variance equal to 0.05 in each coordinate of the ambient space.   (a) The colors represent the true response values defined on the M\"{o}bius Strip. (b) 100 labeled instances (red numbers) and fitted response values (colors on the surface of the M\"{o}bius Strip) by a kernel regression with the D-optimal Design. (c) 100 labeled instances (red numbers) and fitted response values  (colors on the surface of the M\"{o}bius Strip)  by ODOEM. Similarly as before, the fitted function in (c) approximates the true function on the M\"{o}bius Strip in (a) better than the fitted function in (b), with or without noise.}
\end{figure}

\begin{figure}[H]
        \centering
        \begin{subfigure}[b]{0.31\textwidth}
            \centering
            \includegraphics[width=\textwidth]{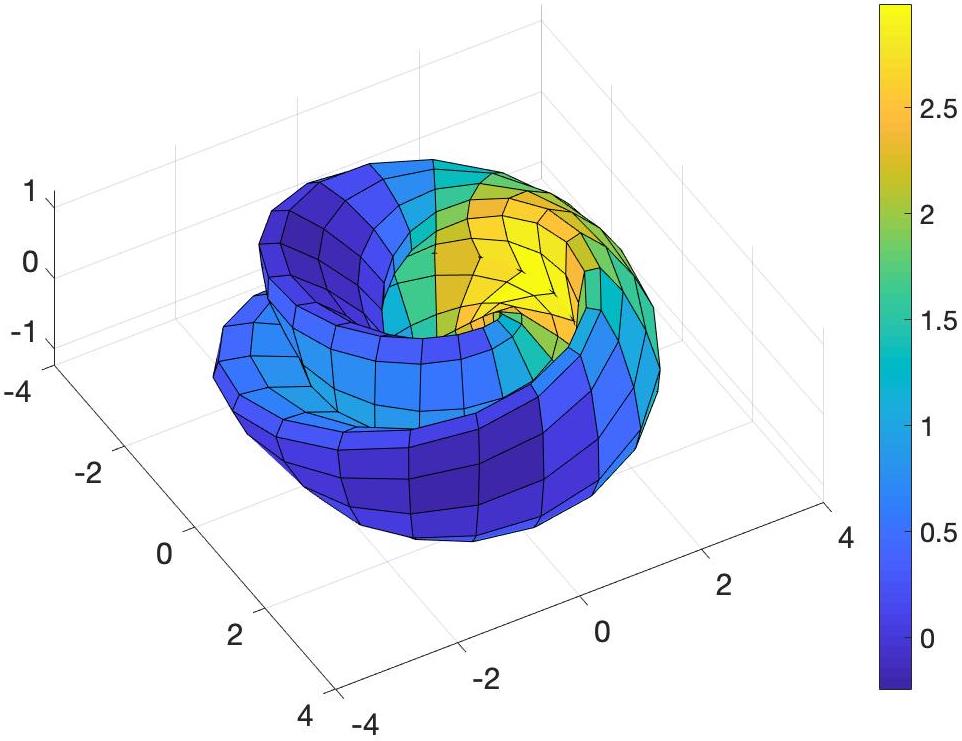}
        \end{subfigure}
        \quad
        \begin{subfigure}[b]{0.31\textwidth}
            \centering
            \includegraphics[width=\textwidth]{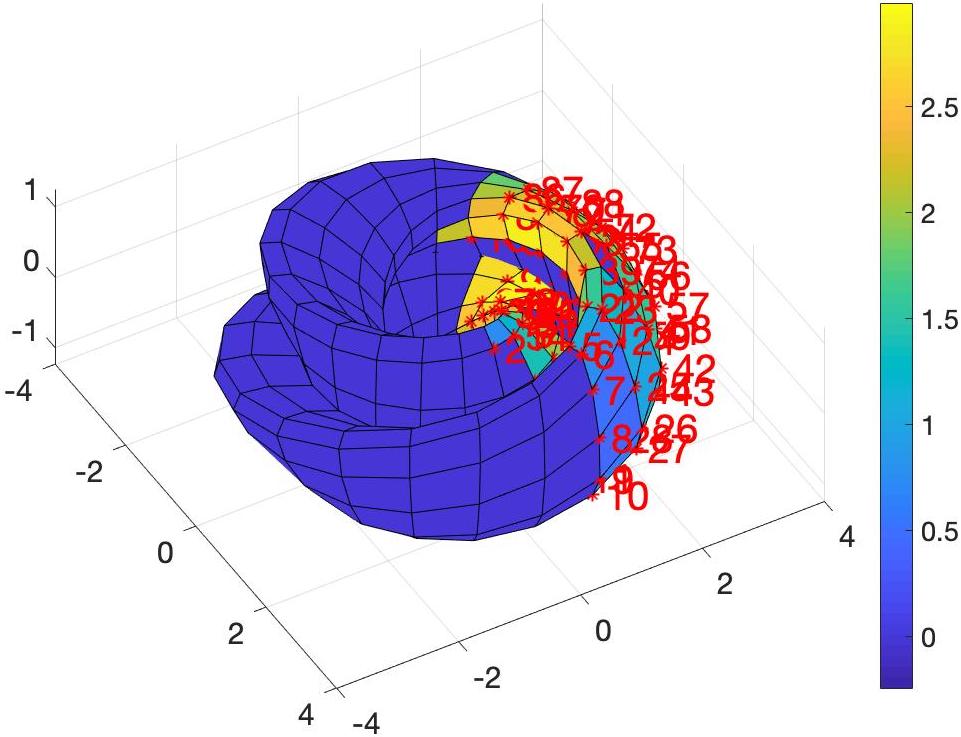}
        \end{subfigure}
        \quad
        \begin{subfigure}[b]{0.31\textwidth}
            \centering
            \includegraphics[width=\textwidth]{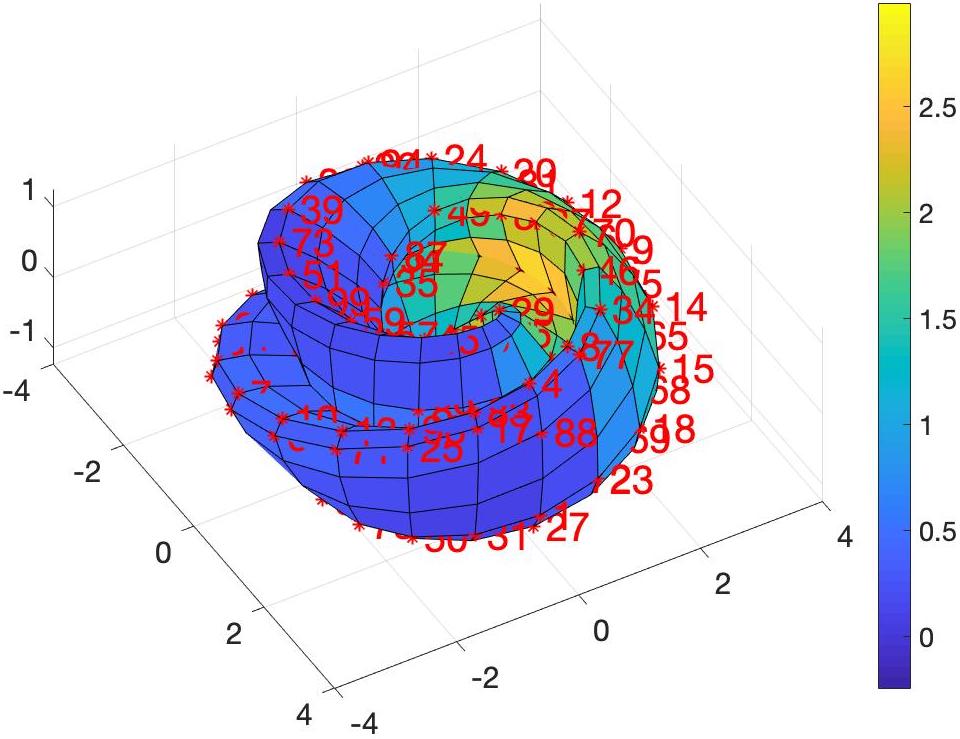}
        \end{subfigure}
         \vskip\baselineskip
         \begin{subfigure}[b]{0.31\textwidth}
            \centering
            \includegraphics[width=\textwidth]{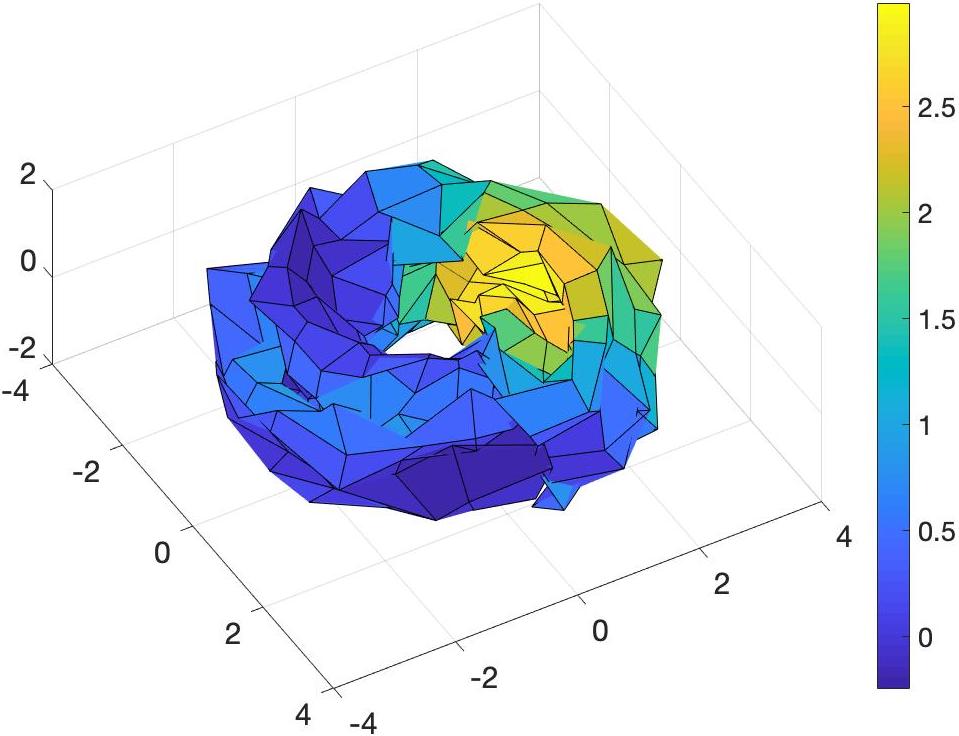}
            \caption{True Function}
        \end{subfigure}
        \quad
        \begin{subfigure}[b]{0.31\textwidth}
            \centering
            \includegraphics[width=\textwidth]{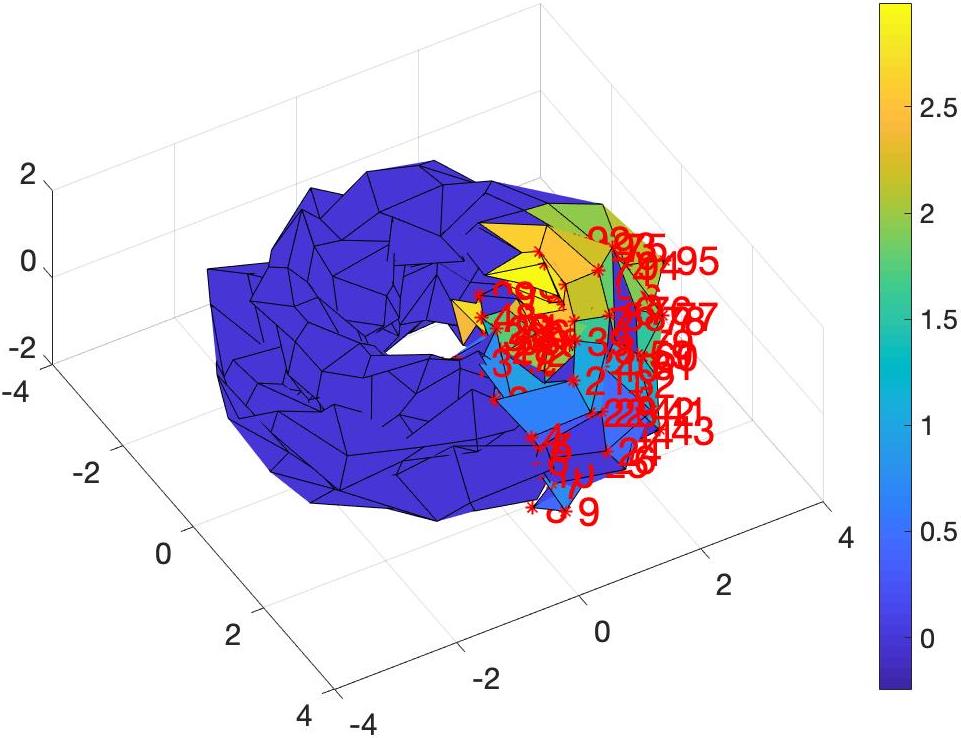}
            \caption{Classical D-optimal}
        \end{subfigure}
        \quad
        \begin{subfigure}[b]{0.31\textwidth}
            \centering
            \includegraphics[width=\textwidth]{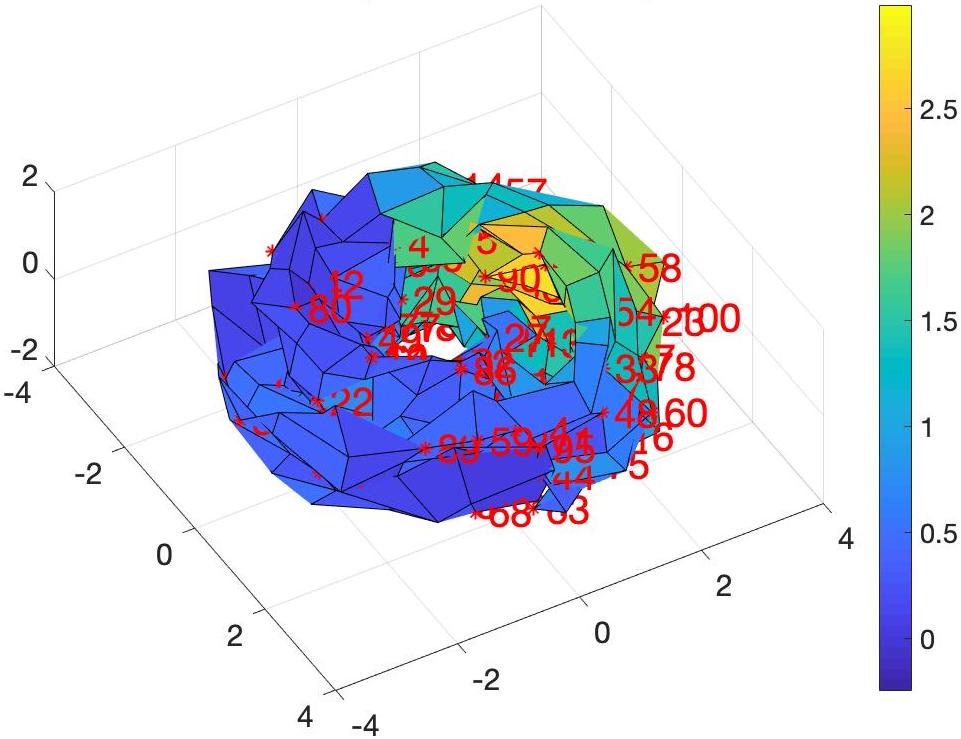}
            \caption{ODOEM}
        \end{subfigure}
        \caption{A figure ``8'' immersion of a Klein Bottle. Top: when $\{\bx_i\}_{i=1}^n$  lie on the figure 8  Immersion. Bottom: when $\{\bx_i\}_{i=1}^n$ are not exactly on the figure 8 immersion due to noise. The isotropic simulated noise follow a normal distribution with zero mean  and variance equal to 0.2 on each coordinate of the ambient space.   (a) The colors represent the true response values defined on the points on the surface of the manifold. (b) 100 labeled instances (red numbers) and fitted response values (colors on the surface of the figure 8 immersion) by a kernel regression with the D-optimal Design. (c) 100 labeled instances (red numbers) and fitted response values  (colors on the surface of the figure 8 immersion) by ODOEM. Once again, the fitted function in (c) approximates the true function on the figure 8 immersion in (a) better than the fitted function in (b), with or without noise.}
\end{figure}

\clearpage
\subsection{Columbia Object Image Library}
To demonstrate application to a more realistic manifold learning problem, we tested the ODOEM algorithm on the Columbia Object Image Library ({COIL-20}). COIL-20 is a database of grey-scale images of 20 different objects and these images were taken at pose intervals of 5 degrees for each object. There are two versions of this database. In this paper, we choose the processed database that contains 1440 normalized images each made of $32 \times 32$ pixels. 

In this set of experiments, the input data $\{\bx_i\}_{i=1}^n$ are the object images and the response values $\{y_i\}_{i=1}^n$ are the corresponding pose angles of these images with respect to the observer. Given an object image, our goal is to estimate the angle of this object in the image. Among 20 different objects, we choose images of four different objects as illustration: a ``Rubber Duck'', a ``Cannon'', a ``Toy Car'' and a ``Piggy Bank''. For each object, we apply the ODOEM algorithm to decide which instances to label and then train the LapRLS model (\ref{eqn:LapRLS}) to predict the angles of the images using the labeled and unlabeled instances. Comparisons were made with the following alternative algorithms: 
\begin{itemize}
    \item Kernel regression model with a classical D-optimal Design;
    \item Kernel regression model with a random sampling scheme;
    \item Kernel regression model with a $L_2$-discrepancy uniform design \citep{FLS06DMCE};
    \item Kernel regression model with a minimax  uniform design \citep{FLS06DMCE};
    \item Kernel regression model with a maximin uniform design \citep{FLS06DMCE};
    \item SVM model with MAED \citep[Manifold Adaptive Experimental Design,][]{CH12IEEETKDE};
    \item SVM model with TED  \citep[Transductive Experimental Design,][]{YZXG2008}.
\end{itemize}

\begin{figure}[!htbp]
   \centering
   \begin{subfigure}[b]{0.47\textwidth}
            \centering
            \includegraphics[width=\textwidth]{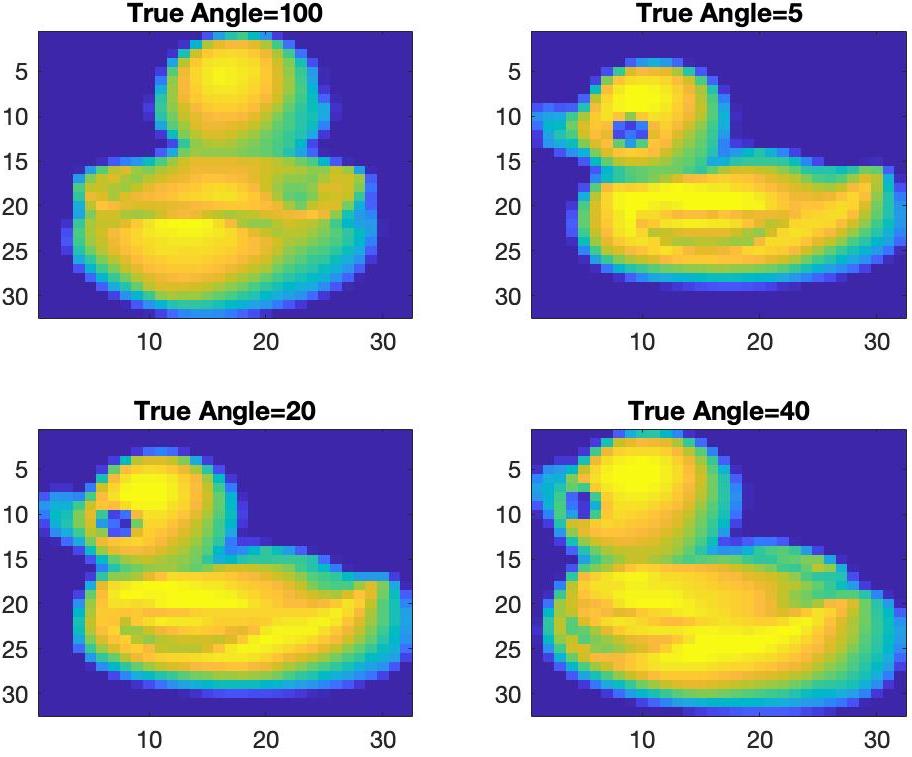}
            \caption{Classical D-optimal Design}
        \end{subfigure}
        \quad
        \begin{subfigure}[b]{0.47\textwidth}
            \centering
            \includegraphics[width=\textwidth]{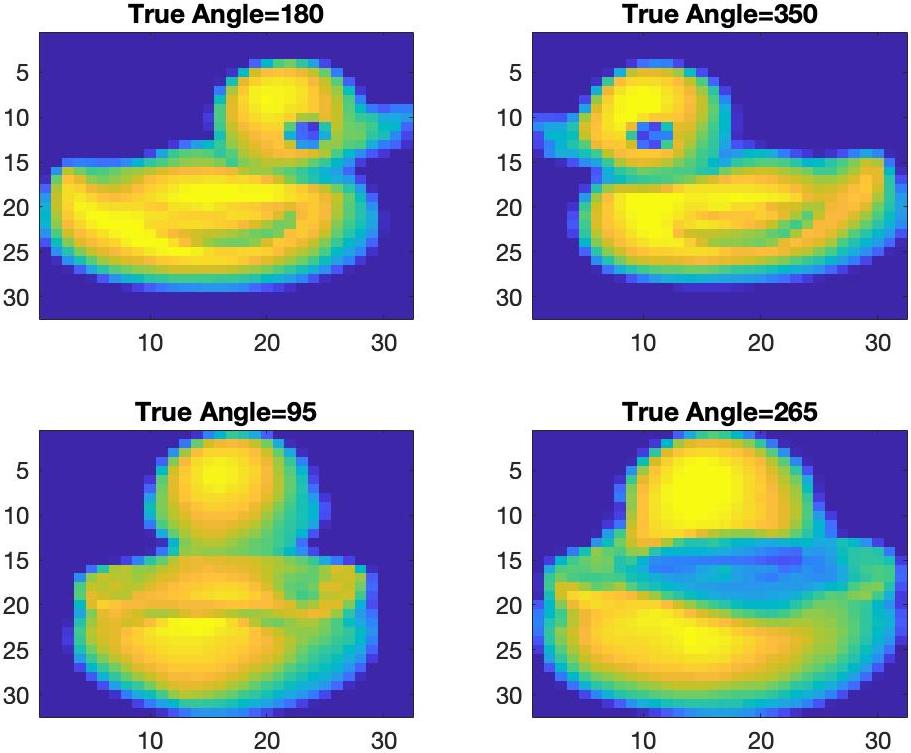}
            \caption{ODOEM}
        \end{subfigure} 
         \vskip\baselineskip
        \begin{subfigure}[b]{0.47\textwidth}
            \centering
            \includegraphics[width=\textwidth]{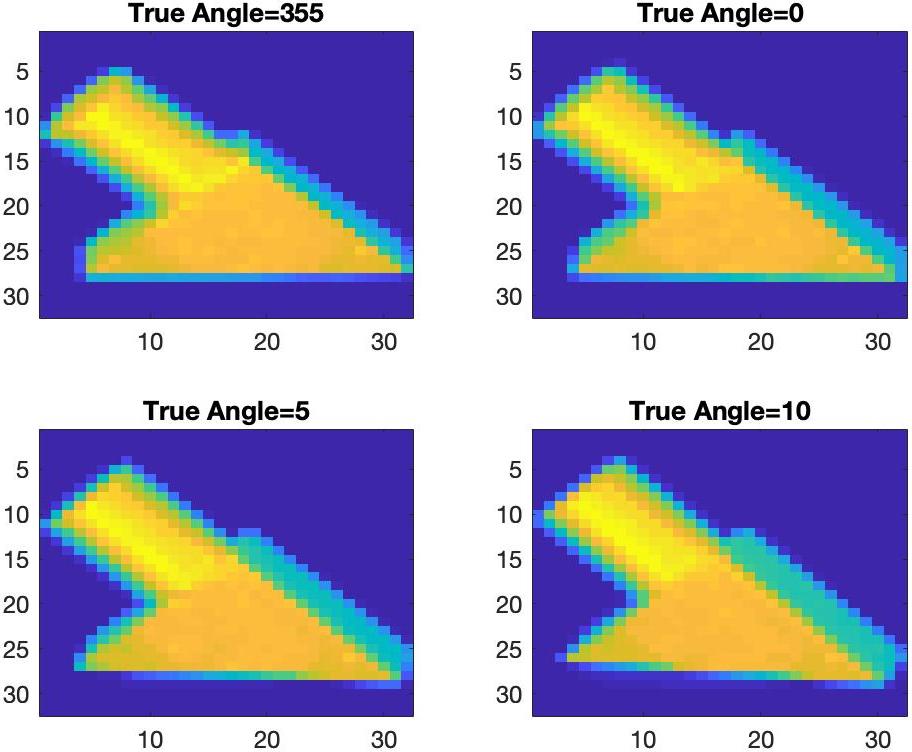}
            \caption{Classical D-optimal Design}
        \end{subfigure}
        \quad
        \begin{subfigure}[b]{0.47\textwidth}
            \centering
            \includegraphics[width=\textwidth]{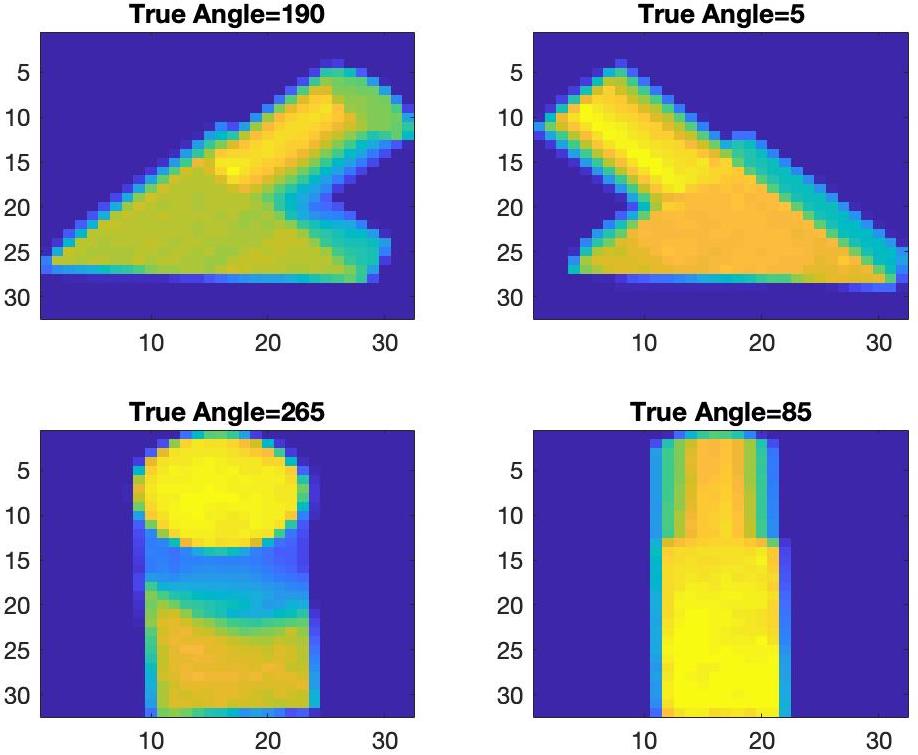}
            \caption{ODOEM}
        \end{subfigure}
    \caption{Top: The first four Rubber Duck images selected by classical D-optimal design and ODOEM. Bottom: The first four Cannon images selected by classical D-optimal design and ODOEM. The true angle is labeled on top of each image. Compared to the classical D-optimal design, there is a greater dispersion (in terms of angles) among the first four images selected by ODOEM.}
    \label{coil20_fig1}
\end{figure}

\begin{figure}[!htbp] 
   \centering
        \begin{subfigure}[b]{0.47\textwidth}
            \centering
            \includegraphics[width=\textwidth]{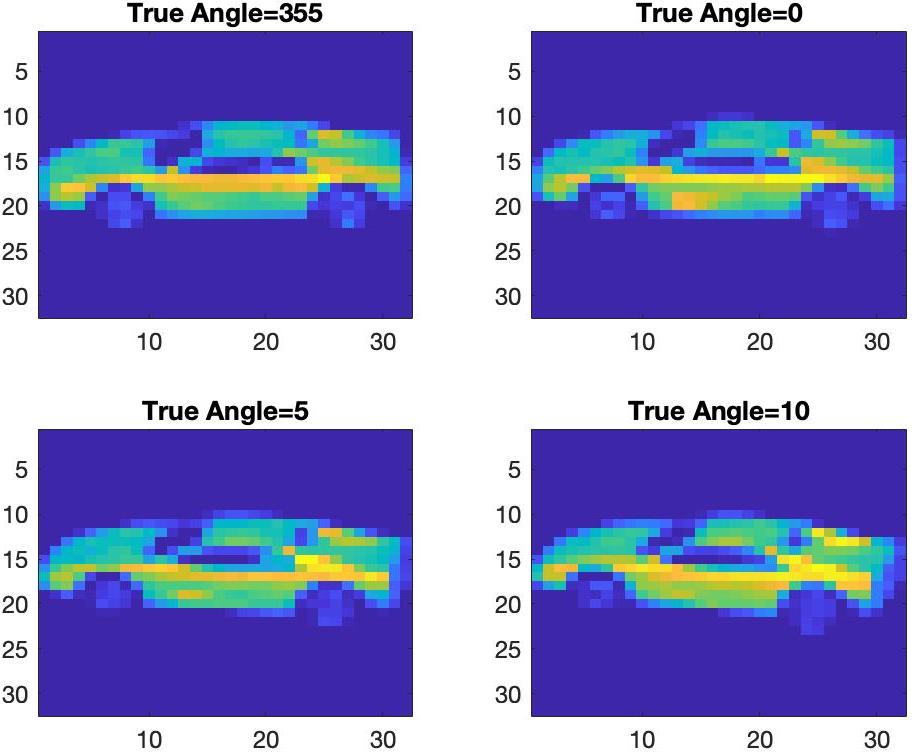}
            \caption{Classical D-optimal Design}
        \end{subfigure}
        \quad
        \begin{subfigure}[b]{0.47\textwidth}
            \centering
            \includegraphics[width=\textwidth]{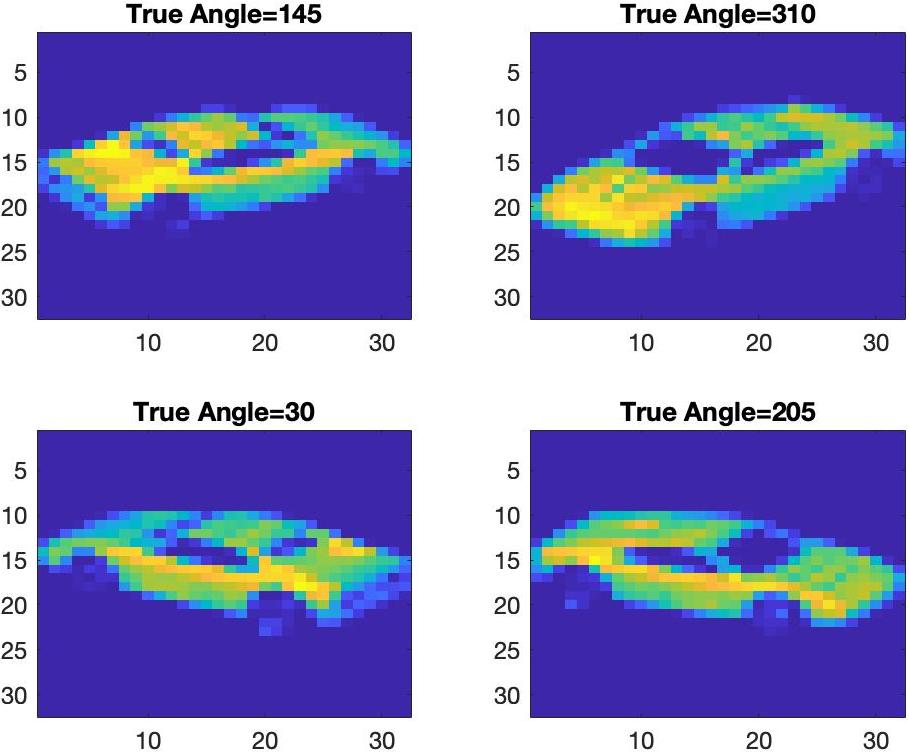}
            \caption{ODOEM}
        \end{subfigure} 
         \vskip\baselineskip
        \begin{subfigure}[b]{0.47\textwidth}
            \centering
            \includegraphics[width=\textwidth]{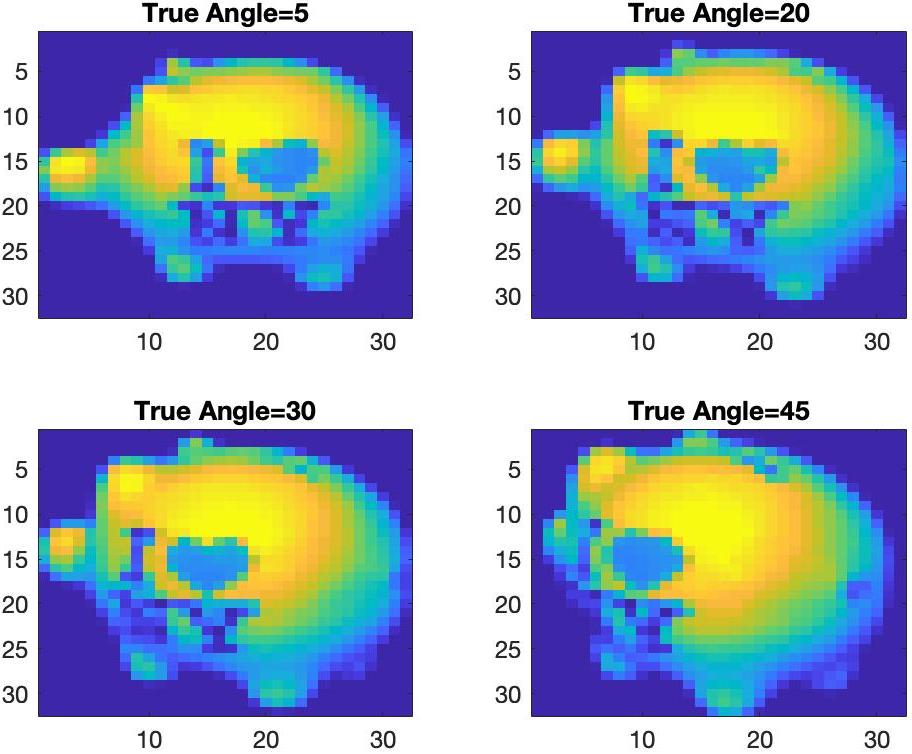}
            \caption{Classical D-optimal Design}
        \end{subfigure}
        \quad
        \begin{subfigure}[b]{0.47\textwidth}
            \centering
            \includegraphics[width=\textwidth]{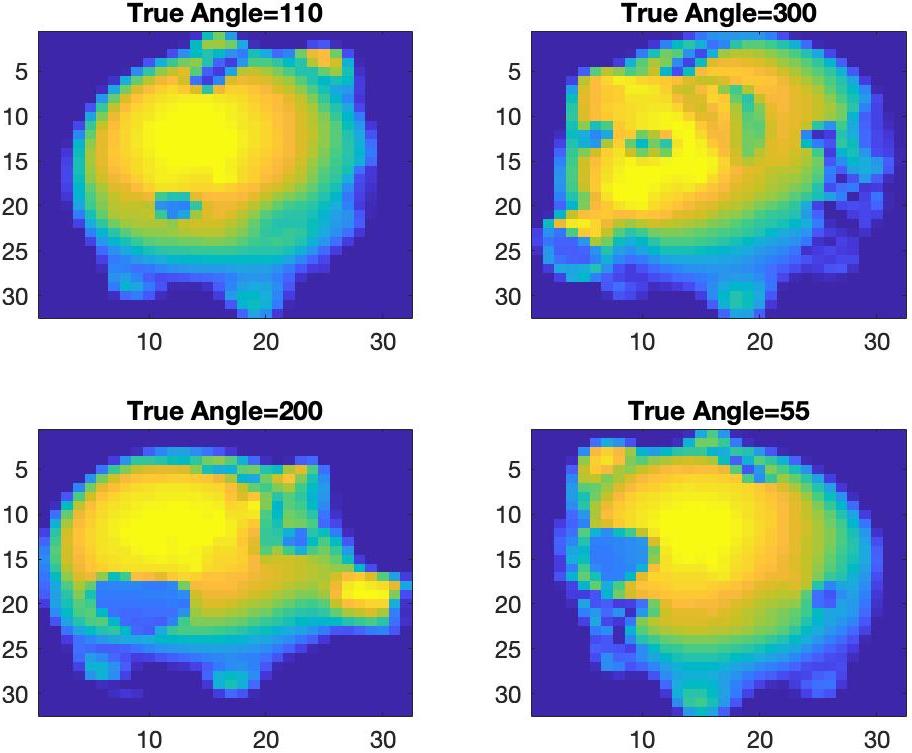}
            \caption{ODOEM}
        \end{subfigure}
    \caption{Top: The first four Toy Car images selected by classical D-optimal design and ODOEM.. Bottom: The first four Piggy Bank images selected by classical D-optimal design and ODOEM. The true angle is labeled on top of each image. Compared to the classical D-optimal design, there is a greater dispersion (in terms of angles) among the first four images selected by ODOEM.}
    \label{coil20_fig2}
\end{figure}

\begin{figure}[!htbp]
        \centering
        \begin{subfigure}[b]{0.475\textwidth}
            \centering
            \includegraphics[width=\textwidth]{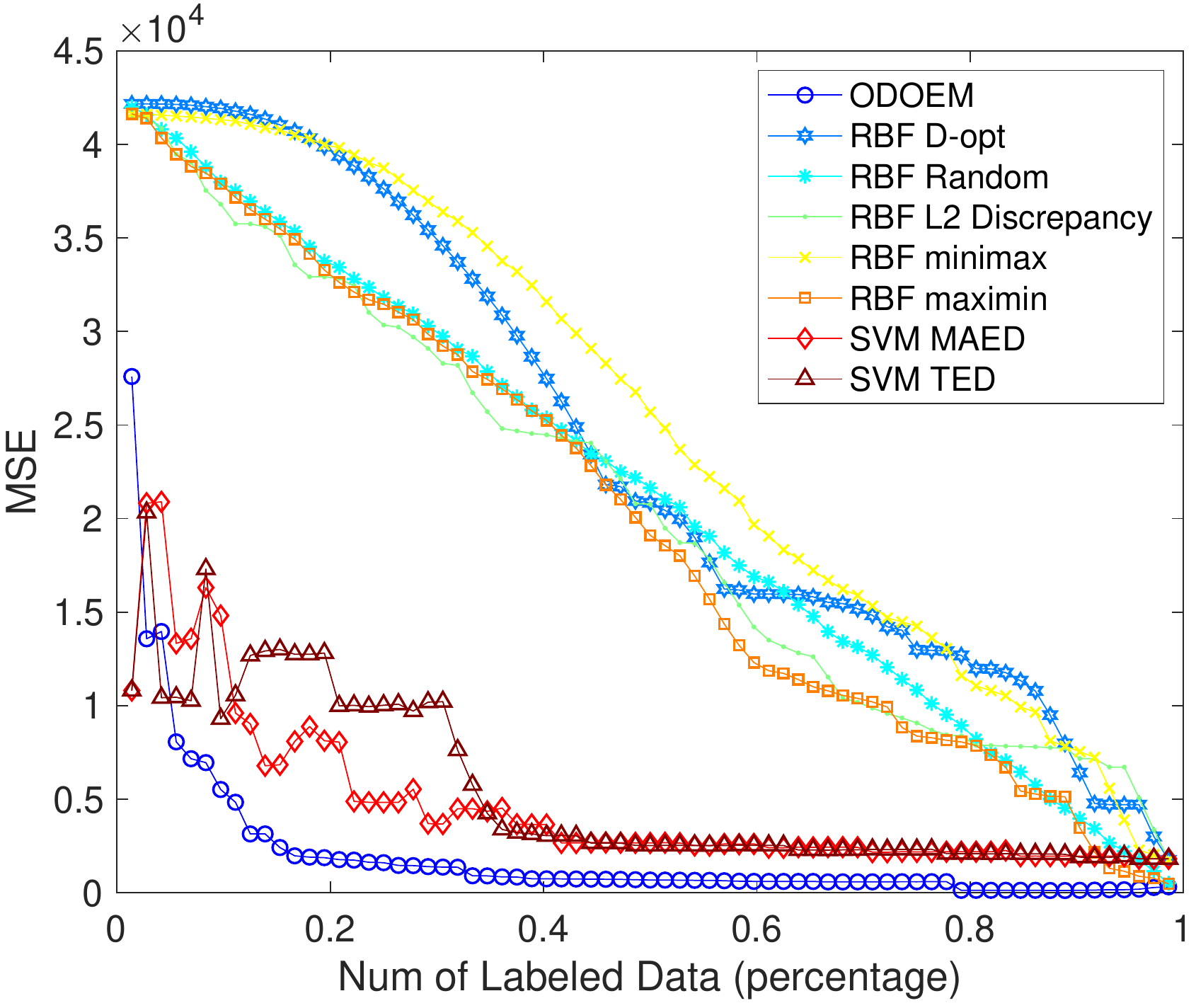}
            \caption{Rubber Duck}
        \end{subfigure}
        \hfill
        \begin{subfigure}[b]{0.475\textwidth}
            \centering
            \includegraphics[width=\textwidth]{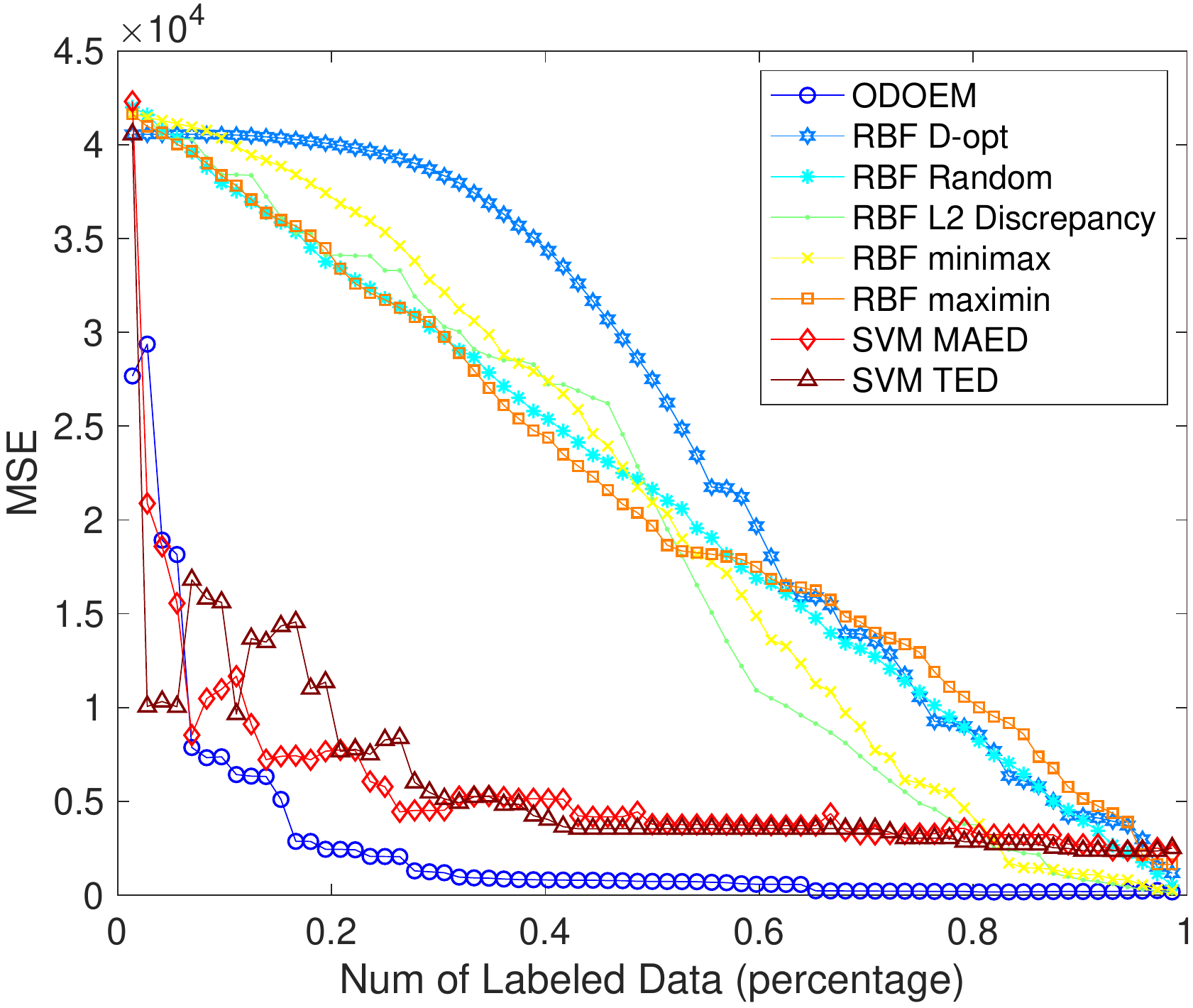}
            \caption{Cannon}
        \end{subfigure}
        \vskip\baselineskip
        \begin{subfigure}[b]{0.475\textwidth}
            \centering
            \includegraphics[width=\textwidth]{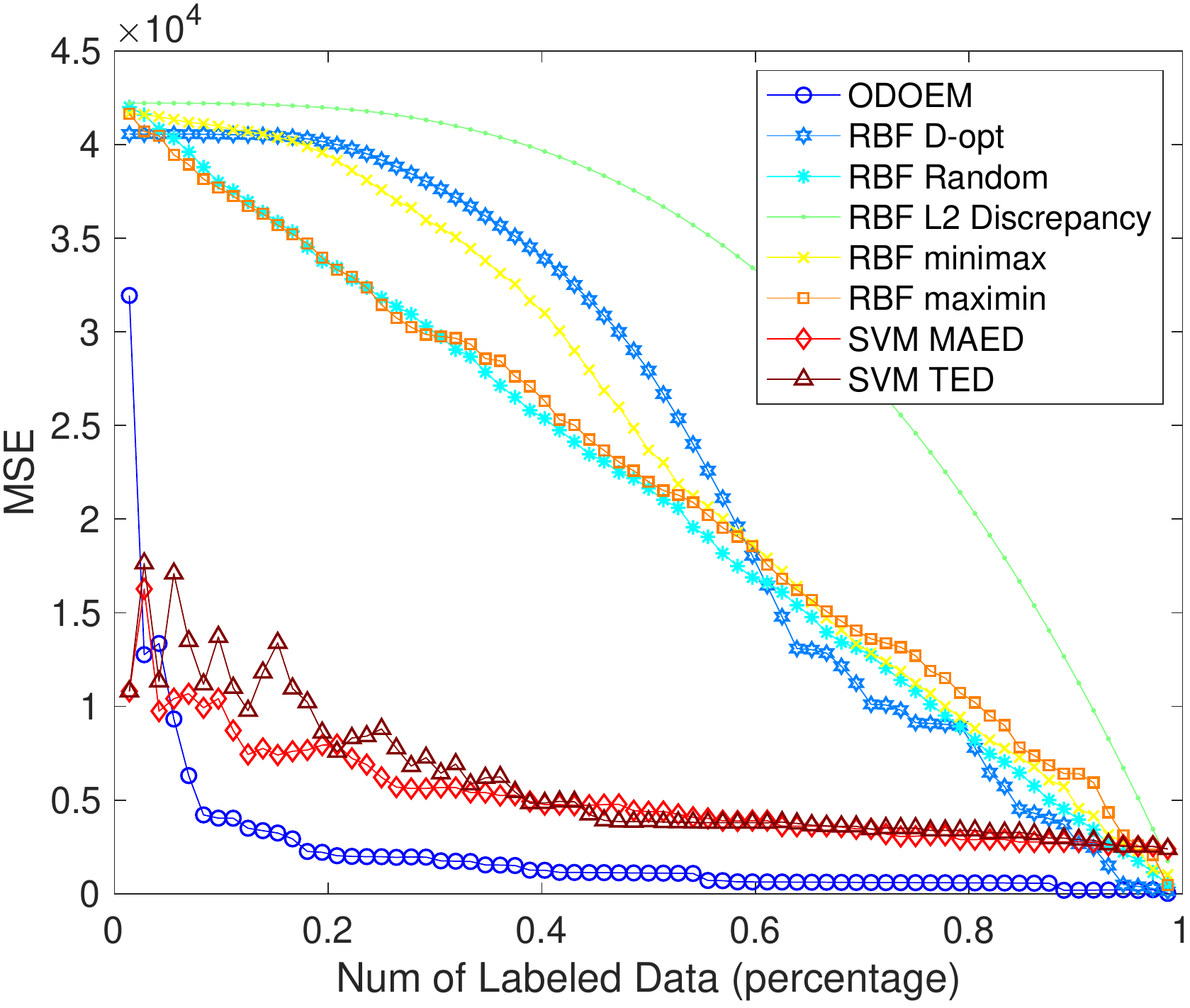}
            \caption{Toy Car}
        \end{subfigure}
        \quad
        \begin{subfigure}[b]{0.475\textwidth}
            \centering
            \includegraphics[width=\textwidth]{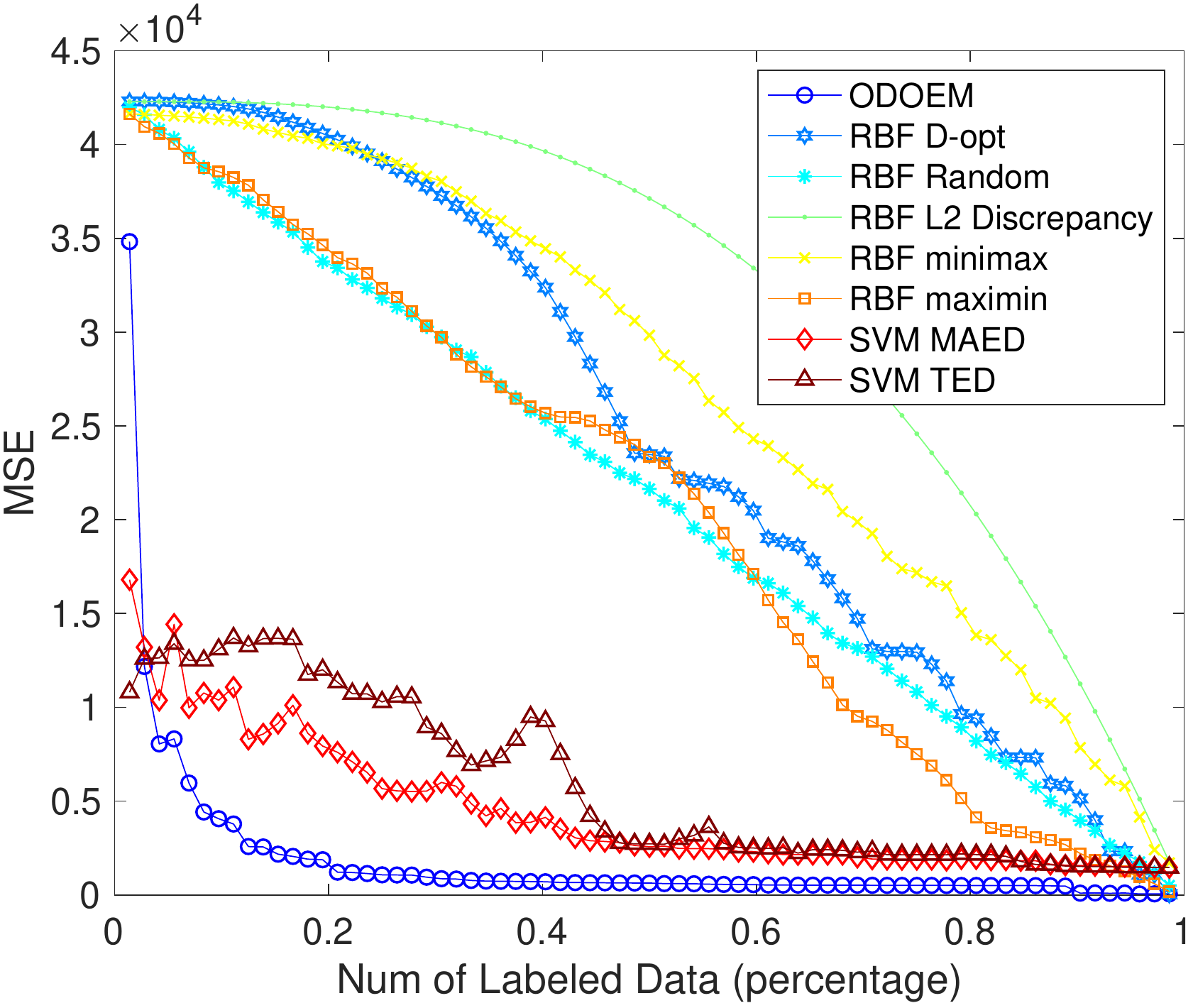}
            \caption{Piggy Bank}
        \end{subfigure}
        \caption{MSE comparison among different algorithms on all four objects. The horizontal axis represents the proportion of images that has been labeled on each object. As it is shown, ODOEM outperforms all the other algorithms on all four object images.}
        \label{coil20_fig3}
\end{figure}

Similar to the synthetic manifold experiments, we used a Radial Basis Function kernel and fixed the range parameter at 0.01 for both kernel regression and SVM. In addition, we choose $\lambda_A=0.01$ and $\lambda_I=-\ln(k/n)$ in ODOEM and kernel regressions. The results are shown in Figures \ref{coil20_fig1}-\ref{coil20_fig3}. In particular, Figures \ref{coil20_fig1} and \ref{coil20_fig2} illustrate the first four images selected by classical D-optimal design and ODOEM for training the models, and Figure \ref{coil20_fig3} demonstrates the fitting performance of different algorithms in terms of MSE.

Based on the results obtained, the following comments can be made: (a) compared to the classical D-optimal design, there is a greater dispersion (in terms of angles) within the first four images selected by ODOEM, which improves the learning curve in Figure \ref{coil20_fig3}; (b) For some uniform design criteria, the corresponding optimization is not convex. Since the images are labeled sequentially, there is no guarantee that the global optimum can be achieved. This explains why some uniform designs do not work very well in these experiments. (c) MAED also benefits from incorporating the manifold structure into the design process. It leads to better fitting performance than most algorithms compared, except ODOEM. (d) ODOEM outperforms all the other algorithms on all four object images.

\section{Conclusions}\label{sec:conc}

In this paper, we have developed a theoretical framework of optimal experimental designs on Riemannian manifolds. Similarly to Euclidean case, we have shown that D-optimal designs and G-optimal designs are also equivalent when the regressors lie on a manifold. Moreover, we have provided a new lower bound for the maximum prediction variance,  demonstrating that this lower bound is achieved at the D/G optimal design. In addition, we proposed a converging algorithm for finding the optimal design of experiments on manifolds. Finally we compared our proposed algorithm with other popular designs and models proposed for both manifold and euclidean optimal design of experiments on several synthetic datasets and real-world image problems, and demonstrated the overall best performance of our ODOEM algorithm.

There are several directions of future research in this work. First, further research can be done to develop a systematic procedure for choosing the regularization parameters $\lambda_A$ and $\lambda_I$. As discussed before, cross-validation is not a feasible strategy in a sequential learning problem since there are few or none labeled instances available at the beginning of the learning process. While we provided heuristic justification for our choices of $\lambda_A$ and $\lambda_I$, a model selection criterion with theoretical guarantees might provide better learning performance. Related work has been discussed by \cite{LCR2019TEST}, where they maximize the likelihood function to choose the values of $\lambda_A$ and $\lambda_I$ in a Gaussian Process model. Furthermore, as it is well-known, there are  optimality criteria other than the D/G criteria in the field of optimal design of experiments. Under different optimal design criteria, new theoretical results of experimental design on manifolds can be explored. Also, for very large scale problems with billions of discrete candidate points, evaluating each point with the corresponding design criteria is computationally exhausting. Some modifications to our algorithm can be investigated, for instance, applying first unsupervised clustering techniques to the covariate data and then evaluate a representative point from each cluster. 

\clearpage
\bibliographystyle{chicago}
\bibliography{DOE}

\begin{thebibliography}{}

\bibitem[\protect\citeauthoryear{Alaeddini, Craft, Meka, and
  Martinez}{Alaeddini et~al.}{2019}]{Alaeddini2019IISE}
Alaeddini, A., E.~Craft, R.~Meka, and S.~Martinez (2019).
\newblock Sequential laplacian regularized v-optimal design of experiments for
  response surface modeling of expensive tests: An application in wind tunnel
  testing.
\newblock {\em IISE Transactions\/}~{\em 51\/}(5), 559--576.

\bibitem[\protect\citeauthoryear{Aronszajn}{Aronszajn}{1950}]{A50TAMS}
Aronszajn, N. (1950).
\newblock Theory of reproducing kernels.
\newblock {\em Transactions of the American Mathematical Society\/}~{\em 68},
  337--404.

\bibitem[\protect\citeauthoryear{Belkin}{Belkin}{2003}]{B03UCTHESIS}
Belkin, M. (2003).
\newblock {\em {Problems of Learning on Manifolds}}.
\newblock Ph.\ D. thesis, The University of Chicago.

\bibitem[\protect\citeauthoryear{Belkin and Niyogi}{Belkin and
  Niyogi}{2003}]{BN03NC}
Belkin, M. and P.~Niyogi (2003).
\newblock Laplacian eigenmaps for dimensionality reduction and data
  representation.
\newblock {\em Neural Computation\/}~{\em 15\/}(6), 1373--1396.

\bibitem[\protect\citeauthoryear{Belkin and Niyogi}{Belkin and
  Niyogi}{2005}]{BN05COLT}
Belkin, M. and P.~Niyogi (2005).
\newblock Towards a theoretical foundation for laplacian-based manifold
  methods.
\newblock In {\em Proceedings of Conference on Learning Theory}.

\bibitem[\protect\citeauthoryear{Belkin, Niyogi, and Sindhwani}{Belkin
  et~al.}{2006}]{BNS06JMLR}
Belkin, M., P.~Niyogi, and V.~Sindhwani (2006).
\newblock Manifold regularization: A geometric framework for learning from
  labeled and unlabeled examples.
\newblock {\em Journal of Machine Learning Research\/}~{\em 7}, 2399--2434.

\bibitem[\protect\citeauthoryear{Box and Draper}{Box and
  Draper}{2007}]{box2007response}
Box, G. and N.~Draper (2007).
\newblock {\em Response Surfaces, Mixtures, and Ridge Analyses}.
\newblock Wiley Series in Probability and Statistics. Wiley.

\bibitem[\protect\citeauthoryear{Cai and He}{Cai and He}{2012}]{CH12IEEETKDE}
Cai, D. and X.~He (2012).
\newblock Manifold adaptive experimental design for text categorization.
\newblock {\em IEEE Transactions on Knowledge and Data Engineering\/}~{\em
  24\/}(4), 707--719.

\bibitem[\protect\citeauthoryear{Chen, Chen, Bu, Wang, Zhang, and Zhang}{Chen
  et~al.}{2010}]{CCBWZZ2010}
Chen, C., Z.~Chen, J.~Bu, C.~Wang, L.~Zhang, and C.~Zhang (2010).
\newblock G-optimal design with laplacian regularization.
\newblock {\em Proceedings of the Twenty-Fourth AAAI Conference on Artificial
  Intelligence\/}~{\em 1}, 413--418.

\bibitem[\protect\citeauthoryear{Cheng and Wu}{Cheng and Wu}{2013}]{CWJASA2013}
Cheng, M. and H.~Wu (2013).
\newblock Local linear regression on manifolds and its geometric
  interpretation.
\newblock {\em Journal of the American Statistical Association\/}~{\em
  108\/}(504), 1421--1434.

\bibitem[\protect\citeauthoryear{Coifman, Lafon, Lee, Maggioni, Nadler, Warner,
  and Zuker}{Coifman et~al.}{2005}]{CLLMNWZ05PNAS}
Coifman, R., S.~Lafon, A.~Lee, M.~Maggioni, B.~Nadler, F.~Warner, and S.~Zuker
  (2005).
\newblock Geometric diffusions as a tool for harmonic analysis and structure
  definition of data: Diffusion maps.
\newblock {\em Proceedings of the National Academy of Sciences\/}~{\em
  102\/}(21), 7426--7431.

\bibitem[\protect\citeauthoryear{Donoho and Grimes}{Donoho and
  Grimes}{2003}]{DG03PNAS}
Donoho, D. and C.~Grimes (2003).
\newblock Hessian eigenmaps: Locally linear embedding techniques for high
  dimensional data.
\newblock {\em Proceedings of the National Academy of Sciences\/}~{\em
  100\/}(10), 5591--5596.

\bibitem[\protect\citeauthoryear{Ettinger, Sangalli, and Perotto}{Ettinger
  et~al.}{2016}]{ESP2016Bio}
Ettinger, B., L.~M. Sangalli, and S.~Perotto (2016, 02).
\newblock {Spatial regression models over two-dimensional manifolds}.
\newblock {\em Biometrika\/}~{\em 103\/}(1), 71--88.

\bibitem[\protect\citeauthoryear{Fang, Li, and Sudjianto}{Fang
  et~al.}{2006}]{FLS06DMCE}
Fang, K., R.~Li, and A.~Sudjianto (2006).
\newblock {\em Design and Modeling for Computer Experiments}.
\newblock Computer Sicence and Data Analysis Series. Chapman and Hall/CRC.

\bibitem[\protect\citeauthoryear{Fedorov}{Fedorov}{1972}]{Fedorov1972}
Fedorov, V.~V. (1972).
\newblock {\em Theory of Optimal Experiments}.
\newblock Academic Press.

\bibitem[\protect\citeauthoryear{Fedorov and Leonov}{Fedorov and
  Leonov}{2013}]{FedorovLeonov2013}
Fedorov, V.~V. and S.~L. Leonov (2013).
\newblock {\em Optimal Design for Nonlinear Response Models}.
\newblock CRC Press.

\bibitem[\protect\citeauthoryear{Gray, Abbena, and Salamon}{Gray
  et~al.}{2006}]{10.5555/1213748}
Gray, A., E.~Abbena, and S.~Salamon (2006).
\newblock {\em Modern Differential Geometry of Curves and Surfaces with
  Mathematica, Third Edition (Studies in Advanced Mathematics)}.
\newblock Chapman \& Hall /CRC.

\bibitem[\protect\citeauthoryear{He}{He}{2010}]{H10IEEETIP}
He, X. (2010).
\newblock Laplacian regularized d-optimal design for active learning and its
  application to image retrieval.
\newblock {\em IEEE Transactions on Imgae Processing\/}~{\em 19\/}(1),
  254--263.

\bibitem[\protect\citeauthoryear{Hein, Audibert, and von Luxburg}{Hein
  et~al.}{2005}]{HAL05COLT}
Hein, M., J.~Y. Audibert, and U.~von Luxburg (2005).
\newblock From graphs to manifolds-weak and strong pointwise consistency of
  graph laplacians.
\newblock In {\em Proceedings of the 18th Conference on Learning Theory}.

\bibitem[\protect\citeauthoryear{Kiefer}{Kiefer}{1974}]{Kiefer1974}
Kiefer, J. (1974).
\newblock General equivalence theory for optimum designs (approximate theory).
\newblock {\em The Annals of Statistics\/}~{\em 2\/}(5), 849--879.

\bibitem[\protect\citeauthoryear{Kiefer and Wolfowitz}{Kiefer and
  Wolfowitz}{1960}]{kiefer_wolfowitz_1960}
Kiefer, J. and J.~Wolfowitz (1960).
\newblock The equivalence of two extremum problems.
\newblock {\em Canadian Journal of Mathematics\/}~{\em 12}, 363–366.

\bibitem[\protect\citeauthoryear{Kimeldorf and Wahba}{Kimeldorf and
  Wahba}{1970}]{kimeldorf1970}
Kimeldorf, G.~S. and G.~Wahba (1970, 04).
\newblock A correspondence between bayesian estimation on stochastic processes
  and smoothing by splines.
\newblock {\em Ann. Math. Statist.\/}~{\em 41\/}(2), 495--502.

\bibitem[\protect\citeauthoryear{Lafon}{Lafon}{2004}]{L04YALETHESIS}
Lafon, S. (2004).
\newblock {\em Diffusion Maps and Geometric Harmonics}.
\newblock Ph.\ D. thesis, Yale University.

\bibitem[\protect\citeauthoryear{Li, Del~Castillo, and Runger}{Li
  et~al.}{2020}]{LCR2019TEST}
Li, H., E.~Del~Castillo, and G.~Runger (2020).
\newblock On active learning methods for manifold data.
\newblock {\em TEST\/}~{\em 29\/}(1), 1--33.

\bibitem[\protect\citeauthoryear{Lin, Thomas, Zhu, and Dunson}{Lin
  et~al.}{2017}]{LTZD2017JASA}
Lin, L., B.~S. Thomas, H.~Zhu, and D.~B. Dunson (2017).
\newblock Extrinsic local regression on manifold-valued data.
\newblock {\em Journal of the American Statistical Association\/}~{\em
  112\/}(519), 1261--1273.

\bibitem[\protect\citeauthoryear{Marzio, Panzera, and Taylor}{Marzio
  et~al.}{2014}]{MPT2014JASA}
Marzio, M.~D., A.~Panzera, and C.~C. Taylor (2014).
\newblock Nonparametric regression for spherical data.
\newblock {\em Journal of the American Statistical Association\/}~{\em
  109\/}(506), 748--763.

\bibitem[\protect\citeauthoryear{Pukelsheim}{Pukelsheim}{2006}]{Pukelsheim2006}
Pukelsheim, F. (2006).
\newblock {\em Optimal Design of Experiments}.
\newblock Society for Industrial and Applied Mathematics.

\bibitem[\protect\citeauthoryear{Roweis and Saul}{Roweis and
  Saul}{2000}]{RS00S}
Roweis, S.~T. and L.~K. Saul (2000).
\newblock Nonlinear dimensionality reduction by locally linear embedding.
\newblock {\em Science\/}~{\em 290}, 2323--2326.

\bibitem[\protect\citeauthoryear{Tenenbaum, de~Silva, and Langford}{Tenenbaum
  et~al.}{2000}]{TSL00S}
Tenenbaum, J.~B., V.~de~Silva, and J.~C. Langford (2000).
\newblock A global geometric framework for nonlinear dimensionality reduction.
\newblock {\em Science\/}~{\em 290}, 2319--2323.

\bibitem[\protect\citeauthoryear{Vuchkov}{Vuchkov}{1977}]{Vuchkov1977}
Vuchkov, I. (1977).
\newblock A ridge-type procedure for design of experiments.
\newblock {\em Biometrika\/}~{\em 64\/}(2), 147--150.

\bibitem[\protect\citeauthoryear{Wahba}{Wahba}{1990}]{W90SM}
Wahba, G. (1990).
\newblock {\em Spline models for observational data}.
\newblock Philadelphia, PA: Society for Industrial and Applied Mathematics.

\bibitem[\protect\citeauthoryear{Wu and Hamada}{Wu and Hamada}{2009}]{Wu2009}
Wu, C.-F.~J. and M.~S. Hamada (2009).
\newblock {\em Experiments: Planning, Analysis, and Parameter Design
  Optimization\/} (2 ed.).
\newblock Wiley.

\bibitem[\protect\citeauthoryear{Yao and Zhang}{Yao and
  Zhang}{2020}]{doi:10.1080/01621459.2019.1610660}
Yao, Z. and Z.~Zhang (2020).
\newblock Principal boundary on riemannian manifolds.
\newblock {\em Journal of the American Statistical Association\/}~{\em
  115\/}(531), 1435--1448.

\bibitem[\protect\citeauthoryear{Yu, Zhu, Xu, and Gong}{Yu
  et~al.}{2008}]{YZXG2008}
Yu, K., S.~Zhu, W.~Xu, and Y.~Gong (2008).
\newblock Non-greedy active learning for text categorization using convex
  transductive experimental design.
\newblock In {\em Proceedings of the 31st annual international ACM SIGIR
  conference on Research and development in information retrieval}, Singapore,
  pp.\  635--642.

\bibitem[\protect\citeauthoryear{Zhu, Liu, Cauley, Rosen, and Rosen}{Zhu
  et~al.}{2018}]{nature2018}
Zhu, B., J.~Z. Liu, S.~F. Cauley, B.~R. Rosen, and M.~S. Rosen (2018).
\newblock Image reconstruction by domain-transform manifold learning.
\newblock {\em Nature\/}~{\em 555\/}(7697), 487--492.

\end{thebibliography}

\clearpage 
\begin{center}
 {\large\bf SUPPLEMENTARY MATERIAL}
\end{center}

\noindent \textbf{Proof of equation (\ref{beta_est2}):}\\
Let $A=(Z_k^\top Z_k + \lambda_A I_p + \lambda_I X^\top LX)^{-1}$. Then:
\begin{eqnarray*}
AX^\top(XAX^\top)^{-1}XZ^\top \mathbf{y}&=&(Z^\top X)^{-1} Z^\top X AX^\top(XAX^\top)^{-1}XZ^\top \mathbf{y}\\
&=& (Z^\top X)^{-1}Z^\top XZ^\top \mathbf{y}\\
&=&Z^\top \mathbf{y}
\end{eqnarray*}
Thus, we have 
\begin{eqnarray*}
X^\top(XAX^\top)^{-1}XZ^\top \mathbf{y}&=& A^{-1} Z^\top \mathbf{y}
\end{eqnarray*}
and therefore equation (\ref{beta_est1}) can be reduced to equation (\ref{beta_est2}). \\
\rightline{$\blacksquare$}

\noindent \textbf{Proof of Proposition 1:} 
\begin{eqnarray*}
 M_{Lap}(\epsilon_{3})&=&\int_{z \in \mathcal{X}} \xi_{3}(z)g(z)g(z)^\top dz +C \\
 &=& \int_{z \in \mathcal{X}} \big [(1-\alpha)\xi_{1}(z)+\alpha \xi_{2}(z)\big] g(z)g(z)^\top dz +C \\
&=& (1-\alpha) \int_{z \in \mathcal{X}} \xi_1(z)g(z)g(z)^\top dz+(1-\alpha) C +\alpha \int_{z \in \mathcal{X}} \xi_{2}(z)g(z)g(z)^\top dz + \alpha C \\
&=& (1-\alpha) M_{Lap} (\epsilon_1) + \alpha M_{Lap} (\epsilon_2) 
\end{eqnarray*}
\rightline{$\blacksquare$}
\clearpage

\noindent \textbf{Proof of Proposition \ref{prop2}:} 

Let $M_{ij}$ be the $(i,j)$ cofactor of the matrix $M_{Lap}(\epsilon_3)$ and let $m_{ij}$ be the $(i,j)$ element of the matrix $M_{Lap}(\epsilon_3)$. Then:

\begin{eqnarray*}
\frac{d \log |M_{Lap}(\epsilon_3)|}{d \alpha} &=&|M_{Lap}(\epsilon_3)|^{-1} \frac{d  |M_{Lap}(\epsilon_3) |}{d \alpha} \\
&=& |M_{Lap}(\epsilon_3)|^{-1} \sum_{i=1}^p \sum_{j=1}^p M_{ij} \frac{d m_{ij}(\alpha)}{d \alpha} \\
&=& \sum_{i=1}^p \sum_{j=1}^p \Big(M^{-1}_{Lap}(\epsilon_3)\Big)_{{ji}} \Big(\frac{d M_{Lap}(\alpha)}{d \alpha}\Big)_{{ij}} \\
&=& \Tr\Big(M^{-1}_{Lap}(\epsilon_3) \frac{d M_{Lap}(\alpha)}{d \alpha} \Big) \\
&=& \Tr\Big\{ M^{-1}_{Lap}(\epsilon_3) \frac{d \big[(1-\alpha) M_{Lap} (\epsilon_1) + \alpha M_{Lap} (\epsilon_2)\big]}{d \alpha} \Big\} \\
&=& \Tr \Big\{ M_{Lap}^{-1}(\epsilon_3)\big [M_{Lap}(\epsilon_2)-M_{Lap}(\epsilon_1)\big] \Big\}
\end{eqnarray*}
\rightline{$\blacksquare$}

\noindent \textbf{Proof of Proposition \ref{prop3}:}

1. 
\begin{eqnarray*}
&& \int_{z\in \mathcal{X}} d (z,\epsilon) \xi(z) dz \\
&=& \int_{z\in \mathcal{X}} g(z)^\top M^{-1}_{Lap}(\epsilon) g(z) \xi(z) dz \\
&=& \int_{z\in \mathcal{X}} \Tr \Big \{g(z)^\top M^{-1}_{Lap}(\epsilon) g(z) \Big \} \xi(z) dz \\
&=& \int_{z\in \mathcal{X}} \Tr \Big \{M^{-1}_{Lap}(\epsilon) \big[g(z) g(z)^\top+C-C \big] \Big \} \xi(z) dz
\end{eqnarray*}
\begin{eqnarray*}
&=& \int_{z\in \mathcal{X}} \Tr \Big \{M^{-1}_{Lap}(\epsilon) \big[g(z) g(z)^\top+C \big ] - M^{-1}_{Lap}(\epsilon)  C  \Big \} \xi(z) dz \\
&=& \int_{z\in \mathcal{X}} \Bigg (\Tr \Big \{M^{-1}_{Lap}(\epsilon) \big[g(z) g(z) ^\top +C \big ] \Big \} - \Tr \Big \{ M^{-1}_{Lap}(\epsilon)  C \Big \}  \Bigg) \xi(z) dz \\
&=& \Tr \Big \{M^{-1}_{Lap}(\epsilon) \Big [\int_{z\in \mathcal{X}} g(z) g(z)^\top\xi(z) dz +C \Big] \Big \} - \Tr \Big \{ M^{-1}_{Lap}(\epsilon)  C \int_{z\in \mathcal{X}} \xi(z) dz\Big \}\\
&=& \Tr \Big \{M^{-1}_{Lap}(\epsilon) \Big [\int_{z\in \mathcal{X}} g(z) g(z)^\top \xi(z) dz +C \Big] \Big \} - \Tr \Big \{ M^{-1}_{Lap}(\epsilon)  C \Big \} \\
&=& \Tr \Big \{M^{-1}_{Lap}(\epsilon) M_{Lap}(\epsilon)\Big \} - \Tr \Big \{ M^{-1}_{Lap}(\epsilon)  C \Big \} \\
&=& p - \Tr \Big \{ M^{-1}_{Lap}(\epsilon)  C \Big \} 
\end{eqnarray*}

2. $\int_{z\in \mathcal{X}} d (z,\epsilon) \xi(z) dz = p - \Tr \Big \{ M^{-1}_{Lap}(\epsilon)  C \Big \} $ implies that $p - \Tr \Big \{ M^{-1}_{Lap}(\epsilon)  C \Big \}$ is the mean value of $d (z,\epsilon)$ for given design $\epsilon$. Thus, we have 
\begin{equation*} 
\max_{z \in \mathcal{X}} d(z,\epsilon) \geq p - \Tr \Big \{ M^{-1}_{Lap}(\epsilon)  C \Big \} 
\end{equation*}
\rightline{$\blacksquare$}

\noindent \textbf{Proof of Proposition \ref{prop4}:} 

Let $\epsilon_1$ and $\epsilon_2$ be two arbitrary designs on the experimental region $\mathcal{X}$ and
let $M_{Lap}(\epsilon_1)$ and $M_{Lap}(\epsilon_1)$ be the corresponding information matrices. Define the set of information matrices on $\mathcal{X}$ as
\begin{equation} 
M_{Lap}(\mathcal{X}):= \{ M_{Lap}(\epsilon)| \xi \in \Xi \}
\end{equation}
where $\Xi$ is the set of all probability measure on $\mathcal{X}$. 
Clearly, $M_{Lap}(\epsilon_1), M_{Lap}(\epsilon_2) \in M_{Lap}(\mathcal{X})$. Based on Proposition \ref{prop1}, we have that
\begin{equation} 
M_{Lap}(\epsilon_3)=(1-\alpha)M_{Lap}(\epsilon_1) + \alpha M_{Lap}(\epsilon_2) \in M_{Lap}(\mathcal{X})
\end{equation}
where $M_{Lap}(\epsilon_3)$ is the information matrix for the design $\epsilon_3=(1-\alpha)\epsilon_1+\alpha \epsilon_2$. This implies that $M_{Lap}(\mathcal{X})$ is a convex set. 

\vspace{0.3cm}
In addition, in order to prove $\log |M_{Lap}(\epsilon)|$ is strictly concave, we also need to show that
\begin{equation} 
\log |(1-\alpha) M_{Lap}(\epsilon_1)+\alpha M_{Lap}(\epsilon_2)| > (1-\alpha) \log |M_{Lap}(\epsilon_1)|+ \alpha \log |M_{Lap}(\epsilon_2)|
\end{equation}
for $\forall \; M_{Lap}(\epsilon_1) \neq M_{Lap}(\epsilon_2)$ and $\forall \alpha \in (0,1)$.
It is known that, for any positive-definite matrices $A$ and $B$, 
\begin{equation} 
|(1-\alpha)A+\alpha B| \geq |A|^{1-\alpha} |B|^\alpha, \; \mathrm{where} \;\alpha \in (0,1),   
\end{equation}
where the equality holds only if $A=B$. Since $M_{Lap}(\epsilon)$ is positive-definite, we have that
\begin{equation} 
|(1-\alpha) M_{Lap}(\epsilon_1)+\alpha M_{Lap}(\epsilon_2)| > |M_{Lap}(\epsilon_1)|^{1-\alpha} |M_{Lap}(\epsilon_2)|^\alpha.
\end{equation}
Therefore, 
\begin{equation*} 
\log |(1-\alpha) M_{Lap}(\epsilon_1)+\alpha M_{Lap}(\epsilon_2)| > (1-\alpha) \log |M_{Lap}(\epsilon_1)|+ \alpha \log |M_{Lap}(\epsilon_2)|
\end{equation*}

\rightline{$\blacksquare$}
 
\noindent \textbf{Proof of Proposition \ref{prop5}:} 

Based on Proposition \ref{prop1}, we have
\begin{eqnarray*} 
M_{Lap}(\epsilon_{k+1})&=&(1-\alpha) M_{Lap}(\epsilon_k)+\alpha M_{Lap}(\epsilon(z)) \\
&=& (1-\alpha) M_{Lap}(\epsilon_k)+\alpha (g(z)g(z)^\top+C) \\
&=& (1-\alpha) \Big [ M_{Lap}(\epsilon_k)+\frac{\alpha}{1-\alpha} g(z)g(z)^\top+ \frac{\alpha}{1-\alpha} C \Big]
\end{eqnarray*}

Then
\begin{eqnarray*} 
|M_{Lap}(\epsilon_{k+1})|&=& (1-\alpha)^p \Big|M_{Lap}(\epsilon_k)(I_p+\frac{\alpha}{1-\alpha}M^{-1}_{Lap}(\epsilon_k)g(z)g(z)^\top+ \frac{\alpha}{1-\alpha} M^{-1}_{Lap}(\epsilon_k) C) \Big|\\
&=& (1-\alpha)^p \Big|M_{Lap}(\epsilon_k) \Big| \Big[1+\frac{\alpha}{1-\alpha} d(z,\epsilon_k)+\frac{\alpha}{1-\alpha} \Tr(M^{-1}_{Lap}(\epsilon_k) C) \Big] 
\end{eqnarray*}


\noindent \textbf{Proof of Proposition \ref{prop6}:} 

Based on Equation (\ref{eqprop4}),  $|M_{Lap}(\epsilon_{k+1})|$ is clearly an increasing function with respect to $d(z,\epsilon_k)$. In order to maximize the value of $ \log |M_{Lap}(\epsilon_{k+1})|$, we choose $z_{k+1}=\operatornamewithlimits{argmax}\limits_{z \in \mathcal{X}} d (z, \epsilon_k)$.
Thus, we have that
\begin{eqnarray*} 
 \log |M_{Lap}(\epsilon_{k+1})|  &=& p \log (1-\alpha)+\log |M_{Lap}(\epsilon_k)| \\
 &&+ \log \Big [1+\frac{\alpha}{1-\alpha} d(z_{k+1},\epsilon_k)+\frac{\alpha}{1-\alpha} \Tr(M^{-1}_{Lap}(\epsilon_k) C) \Big ].
\end{eqnarray*}
It can be shown that
\begin{eqnarray*} 
& &\frac{\partial \log |M_{Lap}(\epsilon_{k+1})| }{\partial \alpha} \\
&=& \frac{d(z_{k+1},\epsilon_k)-(p-\Tr(M^{-1}_{Lap}(\epsilon_k) C))+p\alpha(1-d(z_{k+1},\epsilon_k)-\Tr(M^{-1}_{Lap}(\epsilon_k) C))}{(1-\alpha)[(1-\alpha)+\alpha d(z_{k+1},\epsilon_k)+\alpha \Tr(M^{-1}_{Lap}(\epsilon_k) C)]}.
\end{eqnarray*}
Let $$\frac{\partial \log |M_{Lap}(\epsilon_{k+1})| }{\partial \alpha} \geq 0,$$
then 
\begin{eqnarray*} 
d(z_{k+1},\epsilon_k)-(p-\Tr(M^{-1}_{Lap}(\epsilon_k) C))+p\alpha(1-d(z_{k+1},\epsilon_k)-\Tr(M^{-1}_{Lap}(\epsilon_k) C)) \geq 0.
\end{eqnarray*}
After simplification, we have
\begin{eqnarray} 
\alpha \leq \frac{d(z_{k+1},\epsilon_k)-(p-\Tr(M^{-1}_{Lap}(\epsilon_k) C))}{p[d(z_{k+1},\epsilon_k)-(1-\Tr(M^{-1}_{Lap}(\epsilon_k) C))]}
\end{eqnarray}
Clearly,
\begin{eqnarray} \label{eq:alpha:range}
0 <\alpha \leq \frac{d(z_{k+1},\epsilon_k)-(p-\Tr(M^{-1}_{Lap}(\epsilon_k) C))}{p[d(z_{k+1},\epsilon_k)-(1-\Tr(M^{-1}_{Lap}(\epsilon_k) C))]}
\end{eqnarray}
is the non-decreasing direction for the value of $ \log |M_{Lap}(\epsilon_{k+1})|$. In addition, based on the Proposition \ref{prop3} and Theorem \ref{Thm1}, when the D/G optimal design is not achieved, it is clear that 
\begin{eqnarray}
\frac{d(z_{k+1},\epsilon_k)-(p-\Tr(M^{-1}_{Lap}(\epsilon_k) C))}{p[d(z_{k+1},\epsilon_k)-(1-\Tr(M^{-1}_{Lap}(\epsilon_k) C))]} > 0
\end{eqnarray}
which guarantees the existence of $\alpha$ in Equation (\ref{eq:alpha:range}).

Therefore, $\Big\{|M_{Lap}(\epsilon_k)|\Big\}_k$ is a monotonic increasing sequence. 

\rightline{$\blacksquare$}

 \begin{description}

 \item [MATLAB code:] A GNU zipped tar file containing all the necessary files and code to perform the experiments described in this article. It also includes all datasets used as examples in this article and a simple readme file on how to reproduce the results. 

 \item [COIL-20 data set:] Data set used in the illustration of ODOEM algorithm in Section \ref{sec:numerical}. (COIL20\_angle.mat file)
 
 \end{description}

\end{document}